\documentclass[11pt]{article}      %
\usepackage{amsmath,amssymb, amsthm, eucal,bm,accents,mathtools}  
\usepackage{graphicx, amscd}
\usepackage{framed}

\usepackage{epstopdf}

\usepackage{txfonts}                  
\usepackage[T1]{fontenc}     

\usepackage{scrextend}     

\usepackage{enumerate}       
\usepackage{float}                 
\usepackage{tikz}
\usetikzlibrary{matrix, arrows, calc}
\usetikzlibrary{arrows,shapes,positioning}

\theoremstyle{plain}

\theoremstyle{definition}

\numberwithin{theorem}{section}
\numberwithin{equation}{section}
\numberwithin{figure}{section}
\newcommand{\gaction}[2]{\genfrac{}{}{0.5pt}{}{#1}{#2}%
                        \!\lower2pt\hbox{\rotatebox[origin=c]{-90}{{$\looparrowright$}}}}
\newcommand{\dotfraction}[2]{\genfrac{}{}{0.5pt}{}{#1}{#2}%
                        \!\lower.5pt\hbox{{$\circ$}}}

\usepackage{titlesec}                                                           
\titlelabel{\thetitle.\ }                                                         
\titleformat*{\section}{\fontsize{14pt}{14pt} \bf}        
\usepackage{fancyhdr}
\usepackage{sidecap}  
\usepackage[hang,small,sc]{caption}

\def\-{\hbox{\raisebox{.75pt}{-}}}

\def\mathbi{\boldsymbol}

\usepackage{fancybox}

\def\smalll{\scriptsize}

\def\depth{\delta}

\begin{document}

\title{Apollonian depth and the accidental fractal}

\author{Jerzy Kocik 
    \\ \small Department of Mathematics, Southern Illinois University, Carbondale, IL62901
   \\ \small jkocik{@}siu.edu
}


\date{\small\today}

\maketitle

\begin{abstract}
\noindent
The depth function of three numbers representing curvatures of three mutually  tangent circles is introduced. 
Its 2D plot leads to a partition of the moduli space of the triples of mutually tangent circles/disks
that is unexpectedly a beautiful fractal, the general form of which resembles that of an Apollonian disk packing,
 except that it consists of ellipses instead of circles.
\\[5pt]
\\[5pt]
\end{abstract}

\section{Introduction}


The main purpose of this note is to present a fractal
that results as a partition of the space of tricycles (i.e., three mutually tangent disks)
into regions of constant values of the ``depth function''.
This function measures the depth a tricycle is buried in the Apollonian disk packing that it determines.
The structure of this fractal is unexpected, intriguing, 
and provides a rewarding object for further investigations.

The fractal and most of its properties were discovered with the aid of computer experimentation.
We provide initial observations and a preliminary analysis of the findings.
One of the more intriguing outcomes is the occurrence 
of a deformed version of Farey addition of proper fractions,
Stern-Brocot tree,
and a  property of ellipses analogous to that of Ford circles.

This is an example of visualization where geometry and number theory meet in an interesting way.

~

\section{Basic notions}

Any circle bounds two disks, the {\bf inner} and the {\bf outer}.
The former is contained inside the circle, 
the latter extends outside of it and has infinite area in $\mathbb R^2$.
A disk has curvature $\pm a$ if its bounding circle has radius $1/a$.
The negative curvature is given to the outer disks.
In the following, the disks are {\bf  tangent} if they are {\bf externally tangent}, i.e., 
if they share only one point. 
A {\bf tricycle} is a configuration of three mutually tangent disks (or circles).

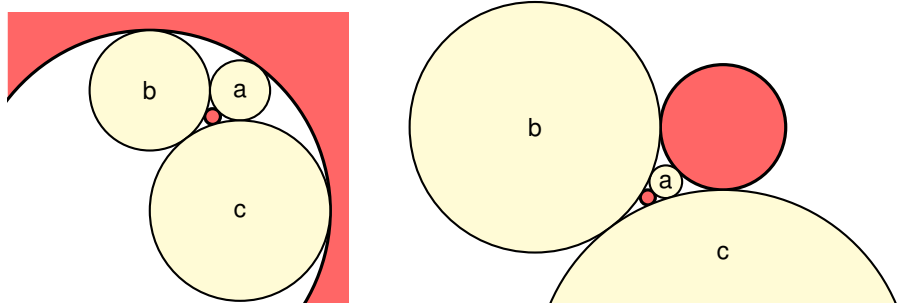
\begin{figure}[h]
\centering
\begin{tikzpicture}[scale=2.4]
\clip (-.79, 1.1) rectangle (1.1, -.51);
\draw [white, fill=red!60] (-.79,-2) rectangle (2,2);
\draw [very thick,  fill=white] (0,0) circle (1);
\draw [very thick,  fill=red!60] (8/23,12/23) circle (1/23);
\draw [thick, fill=yellow!20] (0/3, 2/3) circle (1/3);
\draw [thick, fill=yellow!20] (3/6, 4/6) circle (1/6);
\draw [thick, fill=yellow!20] (1/2, 0/2) circle (1/2);

\node at (3/6, 4/6) [scale=.9, color=black] {\sf a};
\node at (0/3, 2/3) [scale=.9, color=black] {\sf b};
\node at (1/2, 0/2) [scale=.9, color=black] {\sf c};
\end{tikzpicture}
\quad
\begin{tikzpicture}[scale=5, rotate=0]
\clip (-.4,1.05) rectangle (1.01, 0.2);

\draw [very thick, fill=red!60] (15/50, 24/50) circle (1/50);
\draw [very thick,  fill=red!60 ] (3/6, 4/6) circle (1/6);
\draw [thick, fill=yellow!20] (8/23,12/23) circle (1/23);
\draw [thick, fill=yellow!20] (0/3, 2/3) circle (1/3);
\draw [thick, fill=yellow!20] (1/2, 0/2) circle (1/2);

\node at (8/23,12/23) [scale=.9, color=black] {\sf a};
\node at (0/3, 2/3) [scale=.9, color=black] {\sf b};
\node at (1/2, 1/3) [scale=.9, color=black] {\sf c};
\end{tikzpicture}

\caption{Examples of solutions (shown as darker disks) to Descartes' problem 
for disks $a$, $b$, and $c$ .
One of the solutions on the left side has negative curvature.}
\label{fig:twosolutions}
\end{figure}

In a 1643 letter to the princess Elizabeth II of Bohemia,
Ren\'e Descartes proposed 
the following ``{\bf Descartes problem}'':  
given three mutually tangent circles, find the fourth that is simultaneously tangent to all of them.
In the next letter he provided  the solution, known as the {\bf Descartes formula},  
according to which the curvatures $a$, $b$, $c$, $d$ of 
four mutually tangent disks (now called {\bf Descartes configuration}), satisfy:
\begin{equation}
\label{eq:Descartes}
(a+b+c+d)^2 = 2\,(a^2+b^2+c^2+d^2)
\end{equation}
The quadratic nature of Descartes' formula assures in general two solutions:
\begin{equation}
\label{eq:Descartespm}
d = a+b+c \pm 2\sqrt{ab+bc+ca}
\end{equation}
This is consistent with geometry:
there are two different disks that complete a given tricycle to a Descartes configuration, 
as illustrated in Figure \ref{fig:twosolutions}.

\noindent
Another form of Equation \eqref{eq:Descartespm} is a linear relation involving both solutions 
to Descartes problem, say $d$ and $d'$:
\begin{equation}
\label{eq:Descarteslinear}
d +d' \ =\  2a+2b+2c 
\end{equation}

~~

\begin{figure}[H]
\centering
\begin{tikzpicture}[scale=2.2]
\draw (0,0) circle (1);
\foreach \a/\b/\c   in {
1 / 0 / 2 
}
\draw (\a/\c,\b/\c) circle (1/\c)
          (-\a/\c,\b/\c) circle (1/\c);
\foreach \a/\b/\c in {
0 / 2 / 3 ,
0 /4 /15 ,
0 / 6 / 35 , 
0 / 8/ 63,
0 /10 / 99,
0 / 12 / 143
}
\draw[thick] (\a/\c,\b/\c) circle (1/\c)
          (\a/\c,-\b/\c) circle (1/\c) ;
\foreach \a/\b/\c/\d in {
3 / 4 /6, 	8 / 6 / 11,	5 / 12/ 14,	15/ 8 / 18,	8 / 12 / 23,	7 / 24 / 26,
24/	10/	27, 	21/	20/	30, 	16/	30/	35, 	3/	12/	38, 	35/	12/	38, 	24/	20/	39, 	9/	40/	42,
16/	36/	47, 	15/	24/	50, 	 48/	14/	51, 	45/	28/	54, 	24/	30/	59, 	40/	42/	59, 	11/	60/	62,
21/	36/	62, 	48/	28/	63, 	33/	56/	66, 	63/	16/	66, 	8/	24/	71, 	55/	48/	74, 	24/	70/	75,
48/	42/	83, 	80/	18/	83, 	13/	84/	86, 	77/	36/	86, 	24/	76/	87, 	24/	40/	87, 	39/	80/	90
}
\draw (\a/\c,\b/\c) circle (1/\c)       (-\a/\c,\b/\c) circle (1/\c)
          (\a/\c,-\b/\c) circle (1/\c)       (-\a/\c,-\b/\c) circle (1/\c) 
;
\node at (-1/2,0) [scale=1.7, color=black] {\sf 2};
\node at (1/2,0) [scale=1.7, color=black] {\sf 2};
\node at (0,2/3) [scale=1.6, color=black] {\sf 3};
\node at (0,-2/3) [scale=1.6, color=black] {\sf 3};
\node at (1/2,2/3) [scale=1.1, color=black] {\sf 6};
\node at (-1/2,2/3) [scale=1.1, color=black] {\sf 6};
\node at (1/2,-2/3) [scale=1.1, color=black] {\sf 6};
\node at (-1/2,-2/3) [scale=1.1, color=black] {\sf 6};
\node at (8/11,6/11) [scale=.77, color=black] {\sf 1$\!$1};
\node at (-8/11,6/11) [scale=.77, color=black] {\sf 1\!1};
\node at (8/11,-6/11) [scale=.77, color=black] {\sf 1\!1};
\node at (-8/11,-6/11) [scale=.77, color=black] {\sf 1\!1};
\node at (5/14,6/7) [scale=.6, color=black] {\sf 1\!4};
\node at (-5/14,6/7) [scale=.6, color=black] {\sf 1\!4};
\node at (5/14,-6/7) [scale=.6, color=black] {\sf 1\!4};
\node at (-5/14,-6/7) [scale=.6, color=black] {\sf 1\!4};
\end{tikzpicture}
\qquad
%
\begin{tikzpicture}[scale=1.5, rotate=90, shift={(0,2cm)}]  
\clip (-1.25,-2.1) rectangle (1.1,2.2);
\draw (1,-2) -- (1,3);
\draw (-1,-2) -- (-1,3);
\draw (0,0) circle (1);
\draw (0,2) circle (1);
\draw (0,-2) circle (1);
\foreach \a/\b/\c in {
3/4/4,  5/12/12,  7/24/24, 9/40/40  
}
\draw (\a/\c,\b/\c) circle (1/\c)    (\a/\c,-\b/\c) circle (1/\c)
          (-\a/\c,\b/\c) circle (1/\c)    (-\a/\c,-\b/\c) circle (1/\c)   ;
\foreach \a/\b/\c in {
8/    6/   9, 	
15/   8/   16 , 	
24/  20/  25, 	
24/  10/  25, 	
21/  20/  28,
16/	30/	33,    
35/	12/	36,
48/	42/	49,
48/	28/	49,
48/	14/	49,
45/	28/	52,
40/	42/	57,
33/	56/	64,
63/	48/	64,
63/	16/	64,
55/	48/	72,
24/	70/	73,
69/	60/	76,
80/	72/	81,
64/	60/	81,
80/	36/	81,
80/	18/	81
}
\draw (\a/\c, \b/\c) circle (1/\c)          (-\a/\c, \b/\c) circle (1/\c)
          (\a/\c,-\b/\c) circle (1/\c)         (-\a/\c,-\b/\c) circle (1/\c)
          (\a/\c,2-\b/\c) circle (1/\c)       (-\a/\c,2-\b/\c) circle (1/\c)
          (\a/\c,2+\b/\c) circle (1/\c)       (-\a/\c,2+\b/\c) circle (1/\c)
          (\a/\c,-2+\b/\c) circle (1/\c)       (-\a/\c,-2+\b/\c) circle (1/\c)
;
\node at (0,0) [scale=1.9, color=black] {\sf 1};
\node at (0,-2) [scale=1.9, color=black] {\sf 1};
\node at (0,2) [scale=1.9, color=black] {\sf 1};
\node at (-3/4,1) [scale=1.1, color=black] {\sf 4};
\node at (-3/4,-1) [scale=1.1, color=black] {\sf 4};
\node at (3/4,1) [scale=1.1, color=black] {\sf 4};
\node at (3/4,-1) [scale=1.1, color=black] {\sf 4};
\node at (7/8,8/6) [scale=.8, color=black] {\sf 9};
\node at (7/8,-8/6) [scale=.8, color=black] {\sf 9};
\node at (7/8,4/6) [scale=.8, color=black] {\sf 9};
\node at (7/8,-4/6) [scale=.8, color=black] {\sf 9};
\node at (-7/8,8/6) [scale=.8, color=black] {\sf 9};
\node at (-7/8,-8/6) [scale=.8, color=black] {\sf 9};
\node at (-7/8,4/6) [scale=.8, color=black] {\sf 9};
\node at (-7/8,-4/6) [scale=.8, color=black] {\sf 9};
\node at (-0.42,1) [scale=.5, color=black] {\sf 1\!2};
\node at (-0.42,-1) [scale=.5, color=black] {\sf 1\!2};
\node at (.42,1) [scale=.5, color=black] {\sf 1\!2};
\node at (.42,-1) [scale=.5, color=black] {\sf 1\! 2};
\end{tikzpicture}
\caption{Apollonian Window (left) and Apollonian Belt (right)}
\label{fig:Apollo}
\end{figure}
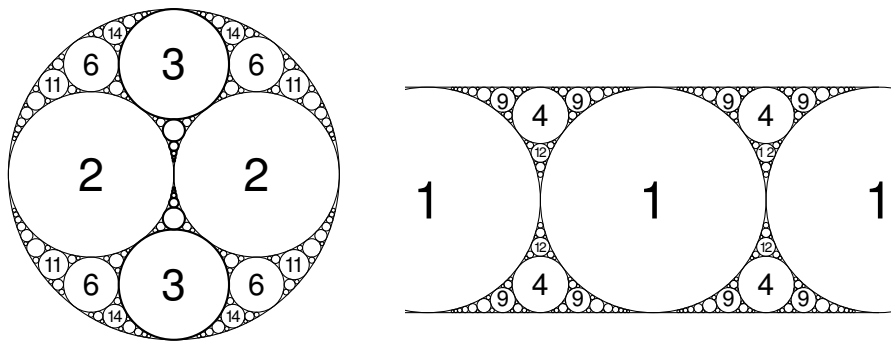

An {\bf Apollonian disk packing} is an arrangement of an infinite number of disks.
Two examples are presented in Figure \ref{fig:Apollo}.
Such an arrangement may be constructed by starting with a tricycle, 
called in this context a {\bf seed},
and completing recursively every tricycle already constructed to a Descartes configuration.

With the help of Descartes' formula (or its linear version), starting with any three tangent circles one may determine 
all curvatures in the resulting Apollonian disk packing.
The disk of a non-positive curvature in a packing will be called the {\bf major disk} 
and its boundary the {\bf major circle}.
\\

\noindent
{\bf Remark:} 
A few clarifications concerning Figure \ref{fig:Apollo} are in order. 
The greatest circle in the Apollonian Window (left) should be viewed as a boundary of an unbounded disk 
extending outwards and having a negative radius and curvature, in this case equal to $(-1)$.  
It is the major disk/circle of of the packing. 
Similarly, the Apollonian Belt is bounded by two lines that should be viewed as the boundaries of half-planes, 
understood as disks of zero curvatures.
This convention allows one to see all pairs of tangent disks in Apollonian packings as tangent {\it externally}.
These two examples of packings  are special: 
they have extra symmetries and 
the curvatures of all disks are integral (see labels inside the disks),
hence the special names mentioned in the figure's caption.
In the following we consider general cases.

\section{Apollonian depth function}

A tricycle determines an Apollonian packing uniquely.
Also, any tricycle in a given Apollonian disk packing
contained in it may serve as its seed.
A natural question arises:
\\\\
{\bf The problem.}
For a given tricycle, how many steps of inscribing new disks are needed to reach the major circle 
of the Apollonian packing that it determines?
Such a number may be viewed as a degree of how deeply a particular triple is buried 
in the network of the packing. 
It is sort of the ``distance'' of the original triple from the external disk
and will be called the ``Apollonian depth'' of the tricycle.
The  goal is to visualize the topological space of tricycles,
and the partition of this space defined by the depth function.
\\

The depth may be found by the following process: 
given a tricycle, form a new one 
by replacing the smallest circle by the greater circle of the two solutions to the Descartes problem \eqref{eq:Descartespm}. 
Repeat until you reach a disk of non-positive curvature.
The number of steps of this process is the value we seek. 
We may give the process an algebraic form without reference to geometry:

~\\
{\bf Definition:}  
The {\bf Apollonian depth} is a function
$$
\depth: \mathbb R^3 \ \to \ \mathbb N \cup \{0,\, \infty\}  
$$
which takes the value zero if any of the three numbers is zero or negative. 
Otherwise, the value is determined by the {\bf dynamical process} in $\mathbb R^3$:
\begin{equation}
\label{eq:dynamical}
T_0 \ \mapsto \  
T_1 \ \mapsto \  
T_2 \ \mapsto \  
T_3 \ \mapsto \  ...
\end{equation}
where $T_0=(a,b,c)$ is the original triple
and each new triple $T_{n+1}$ 
is obtained by replacing the greatest number in $T_n=(a_n,b_n,c_n)$ 
by 
$$
a_n+b_n+c_n-2\sqrt{a_nb_n+b_nc_n+c_na_n}
$$
(The minus sign of  \eqref{eq:Descartespm} to pick the greater disk, see .) \,
The process is to be run until the first occurrence of zero or negative number in some $T_d$.
The number of steps $d$ defines the {\bf Apollonian depth} of the initial  triple $T_0$.
\\

Here is an example:
$$
T_0=(15, 35,102) \ 
              \mapsto\ \underbrace{(15,35, 2)}_{T_1}\ 
              \mapsto\ \underbrace{(15,2,2)}_{T_2}\ 
              \mapsto\ \underbrace{(3,2,2)}_{T_3}\ 
              \mapsto\ \underbrace{(-1,2,2)}_{T_4}
$$
Thus the depth of the triplet $(15,35,102)$ is $\depth(T_0) = 4$.


~\\
{\bf Remark on unbounded packing:}  
Usually, the dynamical system \eqref{eq:dynamical}  terminates after a finite number of steps. 
However, an infinite process is possible!  Consider the following triplet:
\begin{equation}
 \label{eq:golden}
 T_0 \ = \ (\varphi\!-\!\sqrt{\varphi}, \;1, \; \varphi\!+\!\sqrt{\varphi})
\end{equation}
 where $\varphi=\dfrac{1+\sqrt{5}}{2}$ is the golden ratio.
 Denote $p= \varphi-\sqrt{\varphi}$  and $T_0 = (p^{-1}, 1, p )$.
 The reader may check that the dynamical process in this case  takes form:
 $$
(p, 1, p^{-1} ) \ \mapsto \ (p^{2}, p,1 ) 
                       \ \mapsto \ (p^{3}, p^{2}, p ) 
                       \ \mapsto \ (p^{4}, p^{3},p^{2}) \mapsto ...
$$ 
The process never ends with a negative curvature.
In this case we define $\depth(T_0) = \infty$.
For more on this unbounded arrangement of disks consult  \cite{jk-u}.

\section{Visualization of the depth function}
\label{sec:d}

\begin{figure}[H]
\centering
\includegraphics[scale=.6]{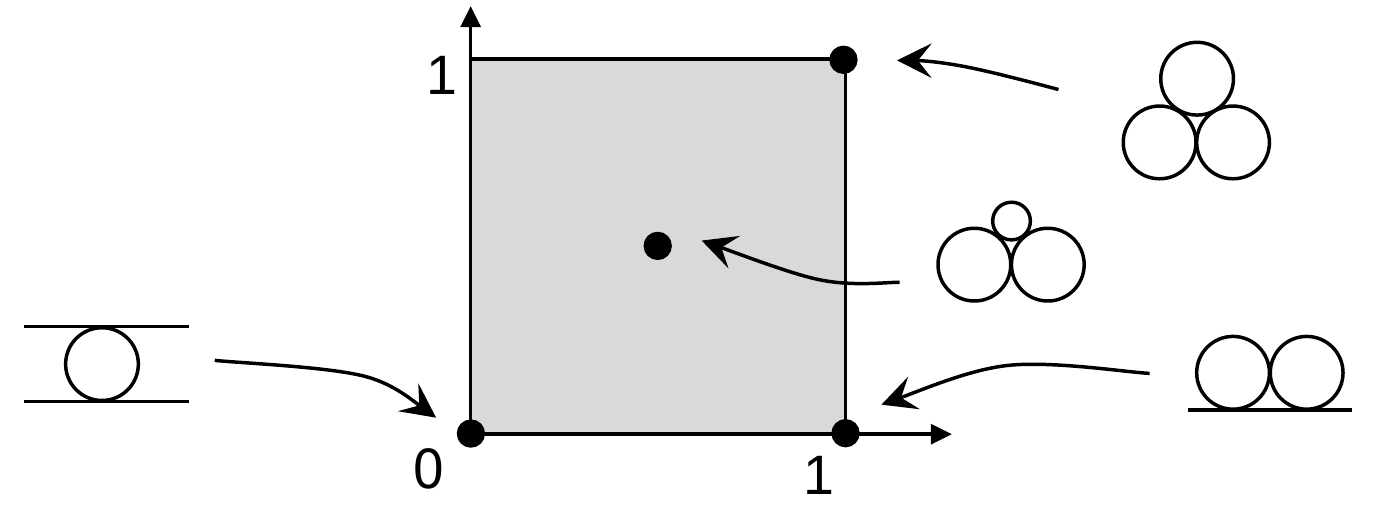} 
\caption{\small Configuration space of tricycles}
\label{fig:configuration}
\end{figure}

The Apollonian depth is invariant under similarity transformations
of the tricycles, i.e., under rotations, translations and dilations.  
In particular:
$$
\depth(a,b,c) = \depth(\lambda a, \lambda b , \lambda c) \qquad \lambda>0
$$
which follows from the homogeneity of Descartes' formula \eqref{eq:Descartes}.
This suggest defining the {\bf moduli space} of tricycles
up to this symmetry group.  As such, it can be 
parametrized by two numbers.
One way is to scale tricycles  
so that its greatest curvature becomes equal to 1.  
Thus the triple $(a,b,c)$ is given coordinates 
\begin{equation}
\label{eq:reduce}
(a,b,c) \ \mapsto \ (x,y) \ = \  \left(\tfrac{a}{c},\tfrac{b}{c}\right)
\end{equation}
where we assumed that $c=\max(a,b,c)$. 
The moduli space $\mathcal T$ of the non-negative triples coincides with the 
unit square $\mathcal T = I^2\subset \mathbb R^2$ parametrized by $(x,y)$
of \eqref{eq:reduce}, see Figure \ref{fig:configuration}.
There is an additional obvious redundancy, $(x,y)\sim (y,x)$,
which will  be ignored for simplicity.

For economy, we shall use the same symbol for this reduced depth function 
$$
\depth (x,y) \ = \  \depth(1, x, y)
$$
The plan is to visualize this function. 
We shall do it by associating to each point of $\mathcal T$  a color or shade 
representing the depth of the corresponding tricycle.



\begin{figure}[h]
\centering
\includegraphics[scale=.34]{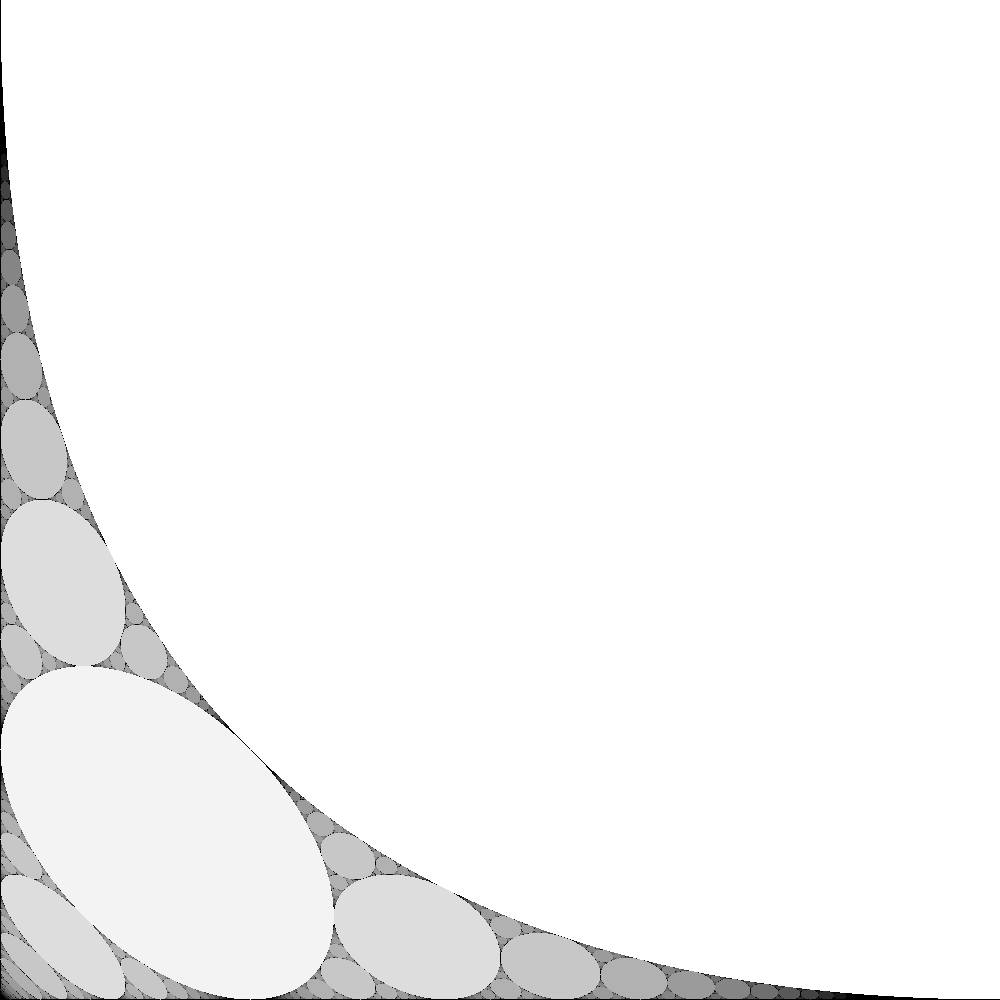}
\caption{\small ``Spiderweb'': Apollonian depth fractal and its unexpected pattern}
\label{fig:chart}
\end{figure}

~\\\\
{\bf Result:}  
The image of this procedure is presented in Figure \ref{fig:chart}
The plot was obtained with the help of the program ``processing.js'' \cite{Processing}.
Appendix A shows the algorithm. 
The computing was done for $1000\times 1000$ points in the square,
uniformly distributed along the rectangular grid. 
The process returned a startling fractal-like pattern which resembles in parts  that of the Apollonian disk packing,
except the disks are replaced by ellipses.
For convenience, we shall refer to it as {\bf spiderweb}, or simply the {\bf web}.
The numbers in Figure \ref{fig:chart-nrs} below  indicate the depth values for selected plateaus.



~

%
%
%

Besides the appealing image of a fractal, more secrets are brought with it,
some of which are explored with the aid of computing.
An alternative representation of the space of tricycles via barycentric coordinates
is presented in Section \ref{sec:barycentric}.


\section{Initial observations}

Here we list some properties 
observed and verified experimentally.

\begin{figure}[H]
\centering
\includegraphics[scale=.35]{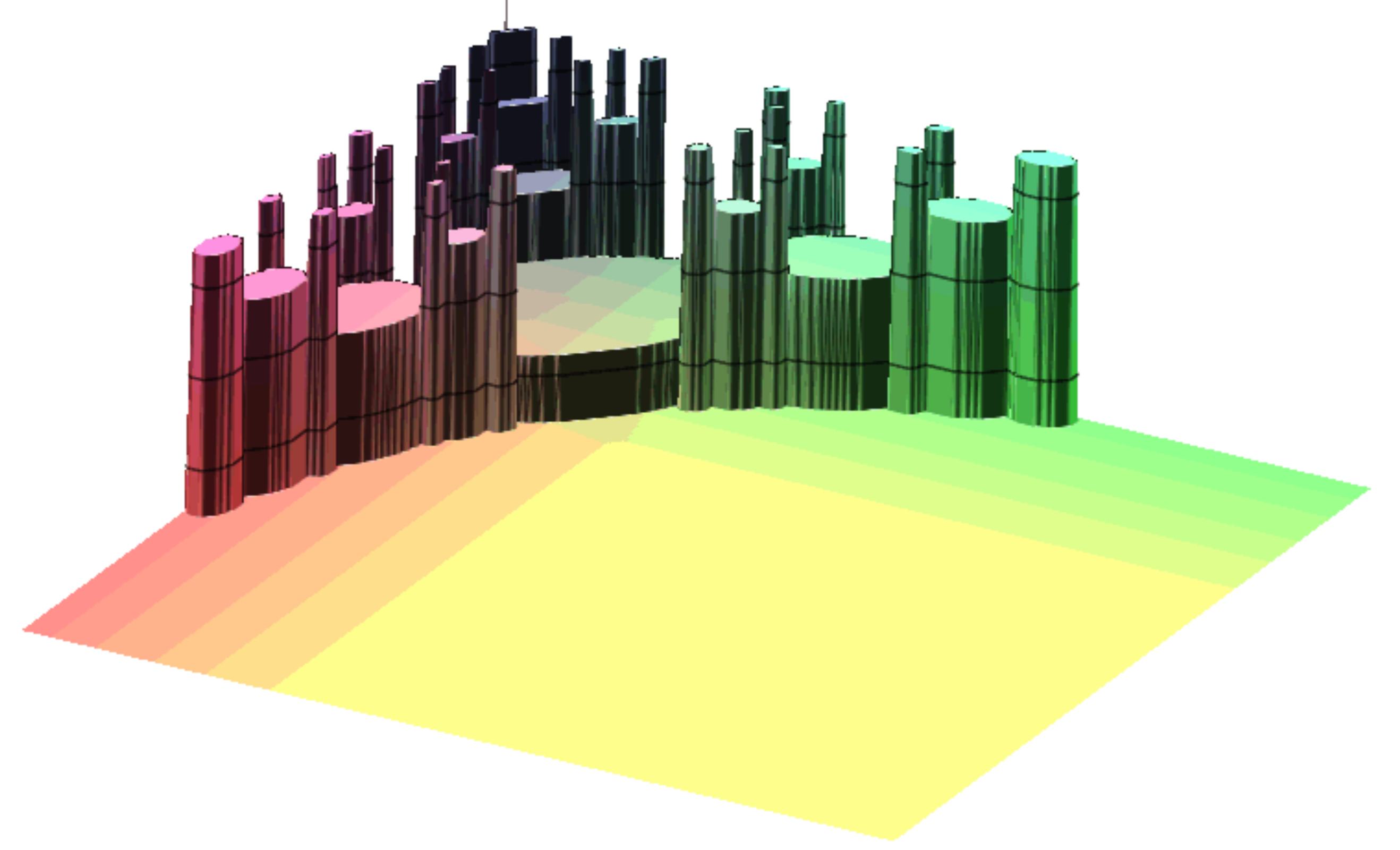}
\caption{The depth function reminds one of a  basalt rock formation 
(figure created with Maple).}
\label{fig:basalt}
\end{figure}

\paragraph{ 1. The general pattern.} 
The spiderweb fractal 
consists of regions composed of points
corresponding to tricycles of the same depth.
The regions turn out to be ellipses (as we will verify) except the the 
main large region in the right upper corner, which is parabolic.
If presented in the standard 3D mode, the graph of the depth function $\delta:\mathbb R^2\to\mathbb R$
has the shape of elliptic columns of different integral heights, quite like geological basalt fields,
see Figure \ref{fig:basalt}. 
Only several regions are included to keep the image clear.
One can see that the plateaus of drastically different heights neighbor each other.

The web may be viewed as a {\bf packing} of the square with ellipses.
It resembles an Apollonian disk packing.  
More accurately --- 
it has the same tangency structure as the circle packing of a square presented in 
Figure \ref{fig:tree}, right.
\footnote{
This is not exactly Apollonian packing since 
not all tangencies follow the Apollonian rule of packing.  
The problem starts with the x-axis and the y-axis being mutually perpendicular, 
and propagates to the ellipses along the diagonal.
However, each region enclosed by two consecutive ellipses 
on the diagonal and the $x$-axis (or $y$-axis) {\it does} follow the Apollonian rule of packing. }

~

\begin{figure}[H]
\centering
\includegraphics[scale=.22]{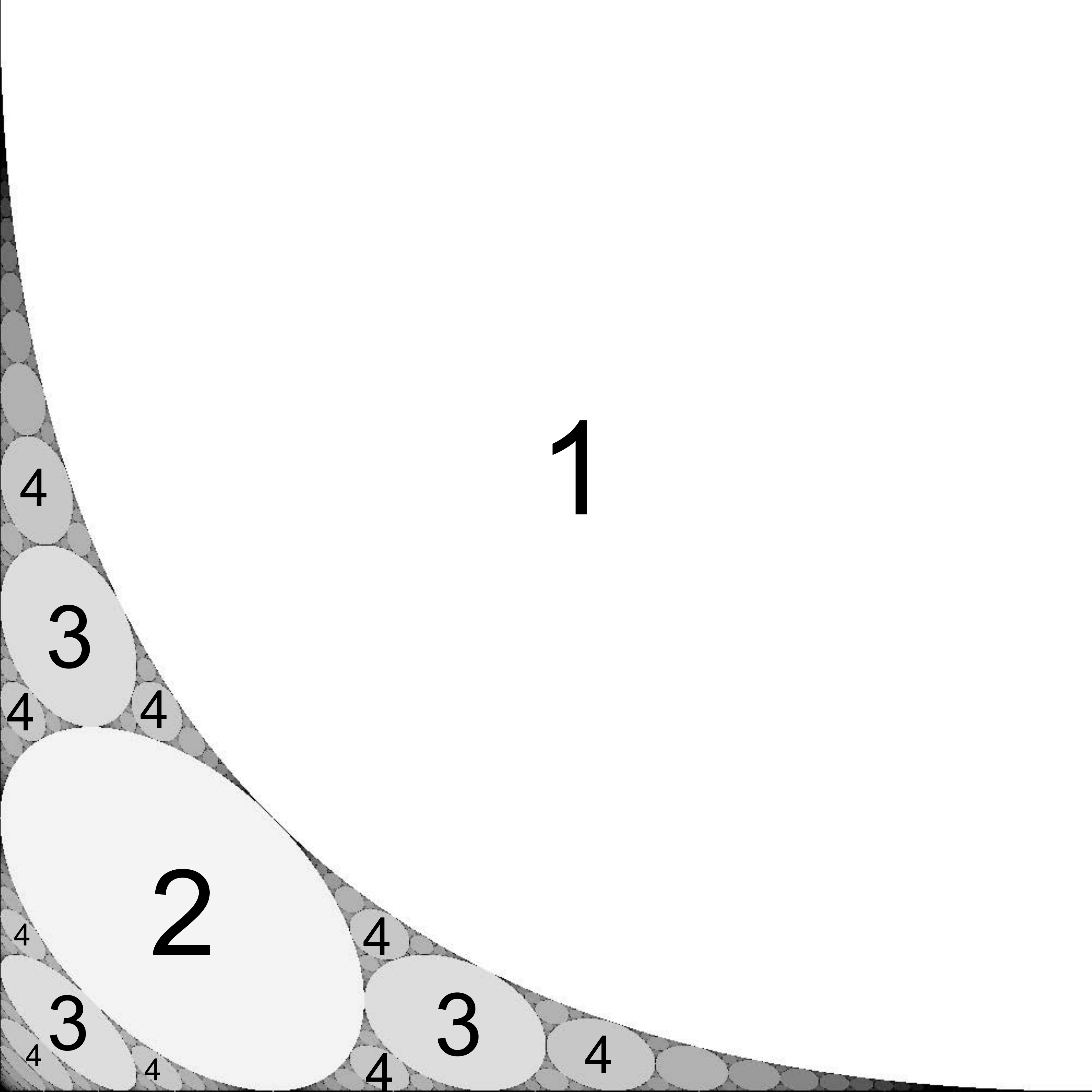}  
\caption{\small Plateaus of constant depth}
\label{fig:chart-nrs}
\end{figure}

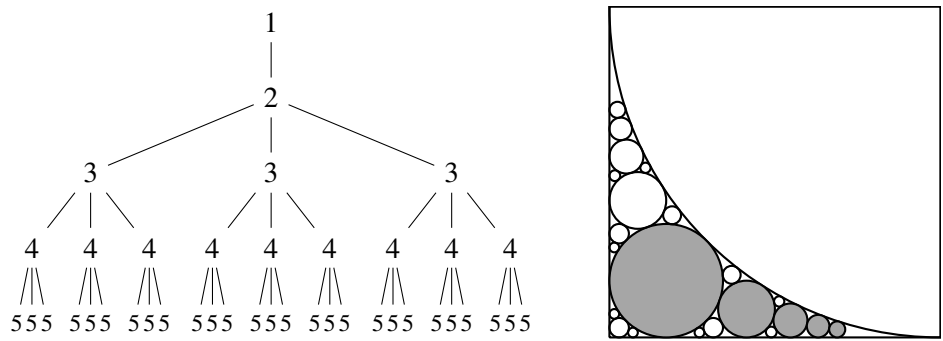
\begin{figure}[H]
\centering
\begin{tikzpicture}[level distance=1cm,
  level 1/.style={sibling distance=1cm},
  level 2/.style={sibling distance=2.4cm},
  level 3/.style={sibling distance=.78cm},
   level 4/.style={sibling distance=.2cm}]
     \node {1}
    child {node {2}
                  child {node {3}
                         child {node {4}
                            child {node{\footnotesize 5}}
                            child{node{\footnotesize 5}} 
                            child {node{\footnotesize 5}} }
                         child {node {4}
                            child {node{\footnotesize 5}} 
                            child {node{\footnotesize 5}}
                            child{node{\footnotesize 5}} }
                         child {node {4}
                           child{node{\footnotesize 5}} 
                           child{node{\footnotesize 5}} 
                           child{node{\footnotesize 5}} }
    }
      child {node {3}
                         child {node {4}
                         child {node{\footnotesize 5}}child{node{\footnotesize 5}} child{node{\footnotesize 5}} }
                           child {node {4}
                           child {node{\footnotesize 5}}child {node{\footnotesize 5}}child{node{\footnotesize 5}} }
                         child {node {4}
                         child{node{\footnotesize 5}} child {node{\footnotesize 5}}child {node{\footnotesize 5}}}
    }
                child {node {3}
                         child {node {4}
                         child {node{\footnotesize 5}}child{node{\footnotesize 5}} child{node{\footnotesize 5}} }
                           child {node {4}
                           child {node{\footnotesize 5}}child{node{\footnotesize 5}} child{node{\footnotesize 5}} }
                          child {node {4}
                          child{node{\footnotesize 5}} child{node{\footnotesize 5}} child{node{\footnotesize 5}} }
      } };
\end{tikzpicture} \quad
\quad
\begin{tikzpicture}[scale=4.4]
\clip (-.01, -.01) rectangle (1.001, 1.001);
\draw [thick] (0,0)--(1,0)--(1,1)--(0,1)--(0,0);  
\draw [thick] (1,1) circle (1); 
\draw [thick, fill=gray!70] (.171573,.171573) circle (.171573); 
\draw [thick] (0.029437, 0.029437) circle (0.029437);  
\draw [thick, fill=gray!70] (0.414213562,0.085786438) circle (0.085786438);
\draw [thick, fill=gray!70] (0.5469181608,  0.05132078829) circle (0.05132078829);
\draw [thick, fill=gray!70] (0.6306019373,  0.03411373215) circle (0.03411373215);
\draw [thick, fill=gray!70] (0.6881924838,  0.02430598181) circle (0.02430598181);
\draw [thick] (0.3137085, 0.0294372) circle (0.0294372);  
\draw [thick] (0.4890416763,  0.01631740053) circle (0.01631740053);
\draw [thick] (0.2720779386,  0.01471862576) circle (0.01471862576);  
\draw [thick] (0.07106781186,  0.01471862576) circle (0.01471862576); 
\draw [thick] (0.3693980626,  0.1893398282) circle (0.02704854689);  
\draw [thick] (0.5128869334,  0.1088913332) circle (0.01555590474); 
%
%
\draw [thick] (0.085786438,0.414213562) circle (0.085786438);
\draw [thick] (0.05132078829, 0.5469181608) circle (0.05132078829);
\draw [thick] (0.03411373215,  0.6306019373) circle (0.03411373215);
\draw [thick] ( 0.02430598181,  0.6881924838) circle (0.02430598181);
\draw [thick] (0.0294372, 0.3137085) circle (0.0294372);
\draw [thick] ( 0.01631740053,  0.4890416763) circle (0.01631740053);
\draw [thick] (0.01471862576,  0.2720779386) circle (0.01471862576);
\draw [thick] (0.01471862576,  0.07106781186) circle (0.01471862576);
\draw [thick] (0.1893398282, 0.3693980626) circle (0.02704854689);
\draw [thick] (0.1088913332,  0.5128869334) circle (0.01555590474);
\end{tikzpicture}
\caption{Left: tree-like structure of the chart.  Right: An analogous circle packing of a square.}
\label{fig:tree1}
\end{figure}

The depth values in elliptic regions is shown in Figure \ref{fig:chart-nrs}.
Note the tree structure of the pattern presented in Figure \ref{fig:tree1} left.
The vertices of this graph correspond to the ellipses,
and the edges join regions that are tangent and differ in depth by~1. 


A few terms will be convenient:
\\[5pt]
$\bullet$  The {\bf main x-wing chain} is the sequence of the ellipses
of depth 2,3,4,..., that are simultaneously tangent to the parabola and 
the $x$-axis.
(The corresponding disks in Figure \ref{fig:tree1} are made dark.)
\\[3pt]
$\bullet$  The {\bf parabolic main wing chain} includes also the ellipses tangent simultaneously to the parabola
and the $y$-axis.
\\[3pt]
$\bullet$  The {\bf diagonal chain}  consists of ellipses simultaneously tangent to 
$x$-axis and the $y$-axis. Their depth values form a sequence 1, 2, 3, 4, 5... .
\\[3pt]
$\bullet$  The {\bf corona} of an ellipse (or parabola) in the web is the set of all 
web ellipses tangent to it.  
\\[3pt]
$\bullet$  The{\bf  $x$-axis corona} consists of all ellipses tangent to the $x$-axis.
The {\bf $y$-axis corona} is defined analogously. 
The {\bf parabolic corona} consists of the ellipse tangent to the parabolic region of depth 1.

\paragraph{2. Warm-up, the first findings.} 
Inspection of 
the $x$-wing chain (regions of depth 2,3,4,... etc)
suggest that
the $x$-coordinate of their points on the $x$ axis seem to follow this simple pattern: 
$$
\frac{1}{3},\quad
\frac{2}{4},\quad
\frac{3}{5},\quad
\frac{4}{6},\ ... \qquad
\hbox{and in general} \ \ \ \ \ 
\frac{n-1}{n+1} \qquad
$$
The guessed formula was verified experimentally with magnification of the plot.
The $y$-values of these points seem to follow also a simple rule: 
$$
\frac{1}{4},\quad
\frac{1}{12},\quad
\frac{1}{24},\quad
\frac{1}{40},\qquad
\hbox{and in general} \
\frac{1}{2n(n+1)}
$$
Thus the series of points separating the ellipses in the main wing are rational:
\begin{equation}
\label{eq:sequence1}
\left(\,  \frac{n-1}{n+1},\,\, \frac{1}{2n(n+1)}\,\right)
\end{equation}

\begin{figure}[h]
\centering
\includegraphics[scale=.37]{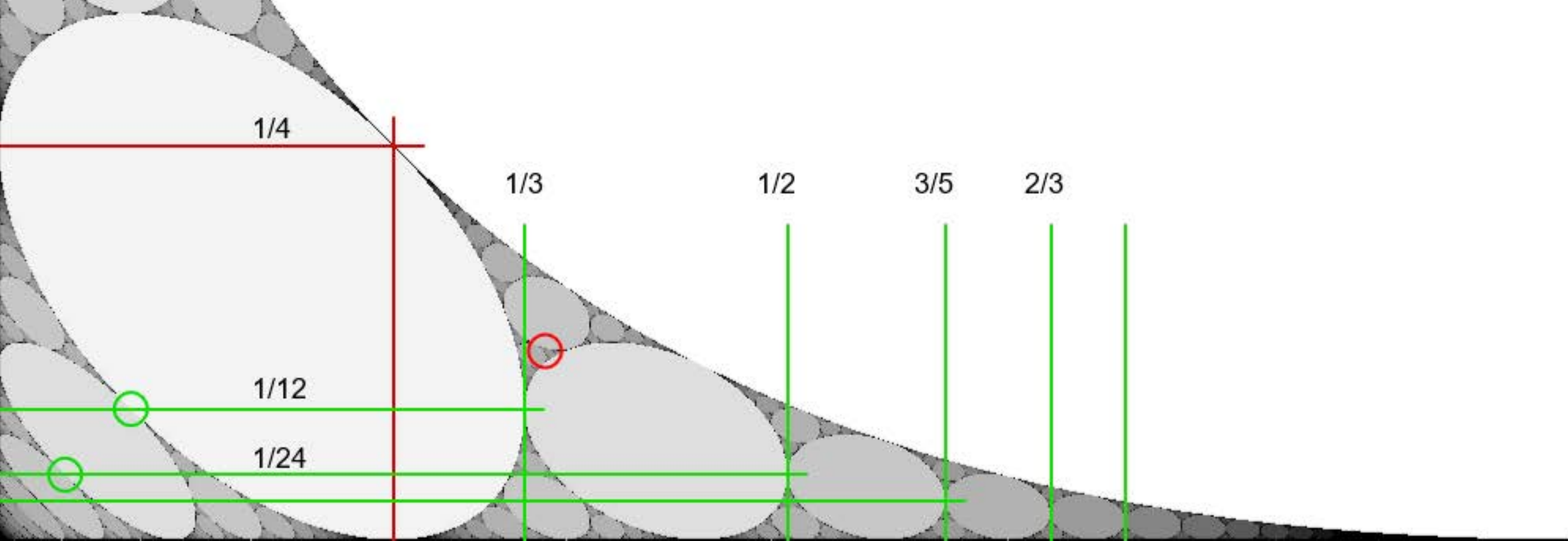}  
\caption{\small Rational coordinates of some tangency points.
Small red circle indicates the golden seed.}
\label{fig:chart-lines}
\end{figure}


An analogous experimental work suggests that 
the sequence of points separating the ellipses along the diagonal of the square $\mathcal T$ are rational 
\begin{equation}
\label{eq:sequence2}
\left(\,  \frac{1}{2n(n+1)},\,\, \frac{1}{2n(n+1)}\,\right)
\end{equation}
Note the common y-values for both sequences, 
\eqref{eq:sequence1} and \eqref{eq:sequence2},
as illuminated by the horizontal lines in Figure \ref{fig:chart-lines}.
All this suggests rationality of the points of tangency in the fractal.
More is coming in Section 2.
\\

\noindent
$\bullet$ 
Yet another puzzling property is visible: each of the ellipses in the central x-chain   
seem to be tangent to the x-axis and the parabola at the same x-coordinate.
\\

\noindent
$\bullet$ 
A property that is harder to notice is that all ellipses in the $x$-axis corona are tangent to the axis 
rational numbers squared.
This will lead to a ``squared'' variation on the Ford fractions. 

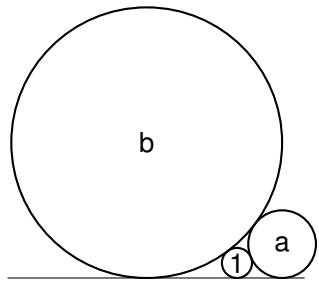
\begin{figure}[H]
\centering
\begin{tikzpicture}[scale=.05]%
\draw [thick] (0, 36) circle (36);
\draw [thick] (2*18, 9) circle (9);
\draw [thick] (2*12, 4) circle (4);
\draw (-37,0) --(42,0);
\node at (0, 36) [scale=1, color=black] {\sf b};
\node at (2*18,9) [scale=1, color=black] {\sf a};
\node at (2*12, 4) [scale=1, color=black] {\sf 1};
\end{tikzpicture}
\caption{Tricycle $(1,\, 4/9,\, 1/9)\sim(9,4,1)$ lies on the boundary between depth 1 and depth 3.
In the chart, it corresponds to the point $(1/9,4/9)$.}
\label{fig:149}
\end{figure}

\paragraph{3. Basalt rock discontinuities}
Intuition would suggest that a sufficiently small change of size of one of the circles 
in a tricycle can lead to a change in depth not greater than 1, if any.
It is not so:  one of the conspicuous outcomes  is that regions 
of arbitrarily big jumps of depth may neighbor each other.
For instance, consider a vertical line through the point separating 
regions of depth $\delta=1$ and region of depth $\delta=3$, namely 
 $(9/16, 1/16)$.
This corresponds to tricycles 
$$
(9/16, \, 1/16 +\varepsilon)
$$
It has depth 1 for $\varepsilon=0$, and 3 for arbitrarily small positive values of $\varepsilon$, 
as is easy to find out with the help of math software. 
The parametrized line $(9/16, t)$ splits into
$$
[0,1] \qquad = \qquad 
\underbrace{\{\,0\,\}}_{\hbox{\footnotesize depth}=0} 
\quad\cup\quad 
\underbrace{(\,0, \,1/9\,)}_{\hbox{\footnotesize depth}=3} 
\quad\cup\quad 
\underbrace{[\,1/9,\, 1\,]}_{\hbox{\footnotesize depth}=1}
$$
Figure \ref{fig:149} illustrates such tricycle. 
Were the disk $b=1/16$ slightly bigger, the line (half-plane) underneath would become a  large disk
and the dynamical system would first acquire a disk to the right of $a$ and $b$ 
before becoming encircleable, hence the depth of the tricycle is 3.

%

%
%
%
%



\begin{figure}[H]
\centering
\begin{tikzpicture}
    \node at (0,0) [anchor=south west,inner sep=0]{\includegraphics[width=\textwidth,height=1.4in]{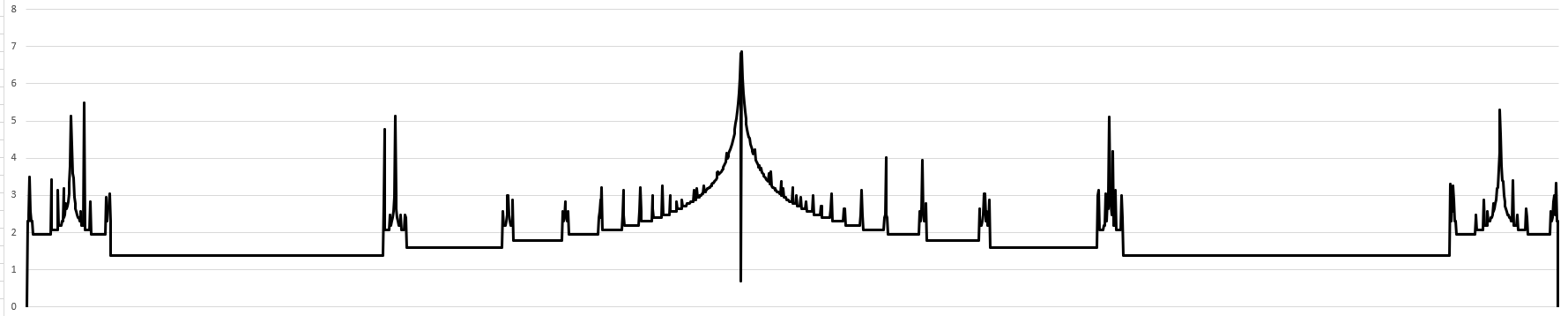}};
    \node at (0.3,-.3) {0};
    \node at (6,-.3) {$\frac{1}{12}$};
    \node at (12.7,-.3) {e};
        \draw (0.28,0.11)--(12.7,0.11);
\end{tikzpicture}
\vspace{-.3in}
\caption{\small A section of the depth function along a vertical segment 
$(\frac{1}{3},0) - (\frac{1}{3},\frac{4-2\sqrt{3}}{3}$).
It consists of 1000 sampling points. 
The vertical axis is in logarithmic scale  Ln($\delta$)
Note the drop at x=1/12 to the value of 2.
 on the vertical axis.}
\label{fig:cut}
\end{figure}
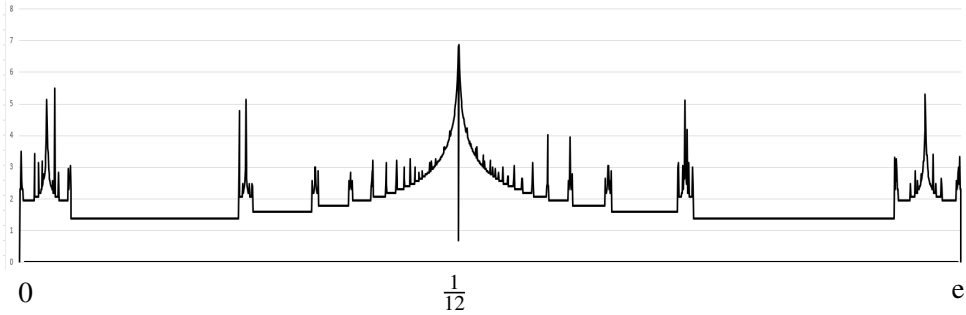

Clearly, one easily finds other points where the tricycle goes through arbitrarily great jumps 
in the value of the depth function $\delta$.
The cylinders of arbitrarily different heights can be mutually tangent in Figure \ref{fig:basalt}.
To exemplify it, a dense sample of a vertical line at $x=1/3$ 
with $y$ changing from 0 to the edge of parabola.
Figure \ref{fig:cut}  presents this cut of the graph of $\delta$,
with logarithmic scale in the value axis.

Here are a few selected points:

$$
\begin{array}{lcl}
\left(1,\, \dfrac{1}{3}, \,\dfrac{1}{12}\right) &\qquad \to \qquad&  \depth= 2\\[11pt]
\left(1,\, \dfrac{1}{3}, \,\dfrac{1}{12}+0.0001\right) &\qquad \to \qquad&  \depth= 836\\[11pt]
\left(1,\, \dfrac{1}{3}, \,\dfrac{1}{12}+0.001\right) &\qquad \to \qquad&  \depth= 86\\[11pt]
\left(1,\, \dfrac{1}{3}, \,\dfrac{1}{12}+0.01\right) &\qquad \to \qquad&  \depth= 11\\[11pt]
\left(1,\, \dfrac{1}{3}, \,\dfrac{1}{12}+0.1\right) &\qquad \to \qquad&  \depth = 1\\
\end{array}
$$

\begin{figure}[H]
\centering
\includegraphics[scale=.12]{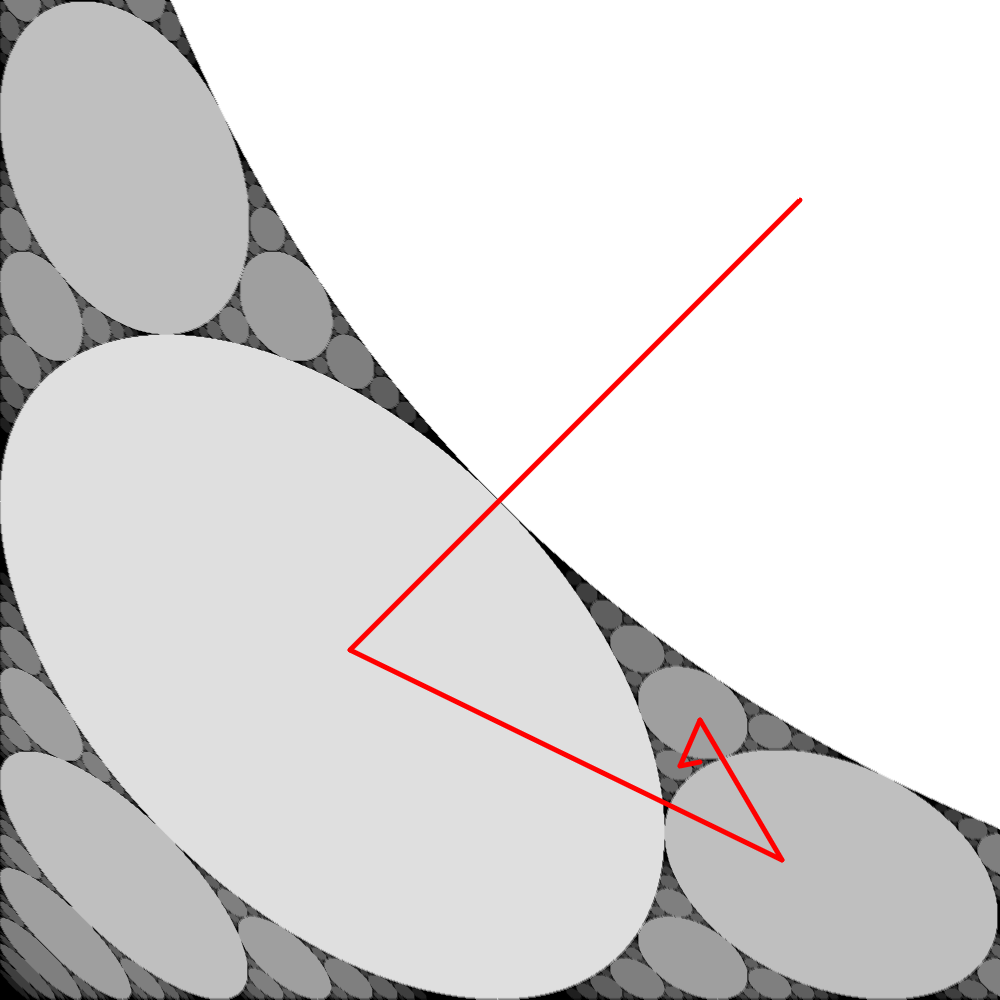}
\caption{\small The spiral that tends to the golden point.
Small red circle indicates the golden seed.}
\label{fig:chart-spiral}
\end{figure}

\paragraph{The golden point (remark):} 
The red circle in Figure \ref{fig:chart-lines}
encircles the point of the unbounded disk arrangement \eqref{eq:golden}, 
for which the depth is equal to infinity.
Its coordinates are:
$$
\left(\, \varphi\!-\!\sqrt{\varphi},\; \left(\varphi\!-\!\sqrt{\varphi}\right)^2\,\right)
\ \approx \ 
\left(\,0.3460,\; 0.1197\,\right)
$$
At first sight the point seems to occupy a generic, unremarkable position in the square.
But a closer look reveals that it is the limit point
of a spiral  that travels through the regions of increasing depth, 1-2-3-...,
each time turning left,
as illustrated in Figure \ref{fig:chart-spiral}.
In Figure \ref{fig:basalt}, these regions form an infinite spiral staircase.


%


Additional experimentation suggests
that points of tangency of two regions of two different depths belong to the
region of the smaller depth.
Thus the regions are in general neither closed nor open.

\section{Calculating the quadratic equations of the plateaus}

Here we derive explicitly the equations for some of the elliptic regions in the web
and show a general method for such calculations.

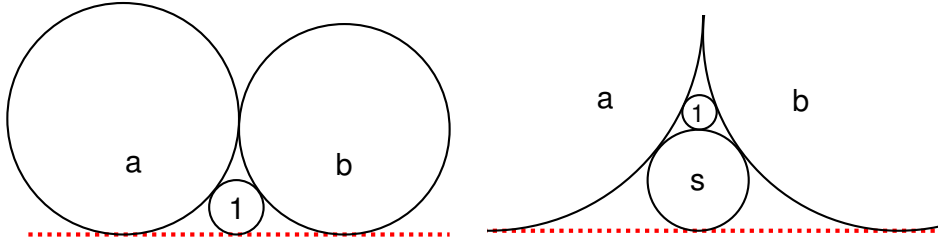
\begin{figure}[h]
\centering
\begin{tikzpicture}[scale=1.4]
\draw [ultra thick, dotted, red] (-2,0)--(2,0);
\draw [thick] (-1.1,1.1) circle (1.1);
\draw [thick] (1,1) circle (1);
\draw [thick] (-.025,.2585) circle (.2585);
\node at (-1, 2/3) [scale=1.1, color=black] {\sf a};
\node at (1, 2/3) [scale=1.1, color=black] {\sf b};
\node at (-.02, 1/4) [scale=1.1, color=black] {\sf 1};
\end{tikzpicture}
\quad
\begin{tikzpicture}[scale=2.6]
\clip (-1.1, -.03) rectangle (1.2, 1.1);
\draw [ultra thick, dotted, red] (-2,0)--(2,0);
\draw [thick] (-1.1,1.1) circle (1.1);
\draw [thick] (1,1) circle (1);
\draw [thick] (-.025,.2585) circle (.2585);
\draw [thick] (-.0182,.606) circle (.088);
\node at (-.5, 2/3) [scale=1.1, color=black] {\sf a};
\node at (.5, 2/3) [scale=1.1, color=black] {\sf b};
\node at (-.028, .6) [scale=.9, color=black] {\sf 1};
\node at (-.028, .25) [scale=1.1, color=black] {\sf s};
\end{tikzpicture}
\caption{Defining the shape of the regions.
Left: $\delta=1$, Right: $\delta=2$.}
\label{fig:shape}
\end{figure}

\paragraph{A. The parabolic plateau of depth 1.}
Referring to Figure \ref{fig:shape}(left), the tricycle $(1,a,b)$ is visibly of depth 1,
but is in the state of a tipping point.
Indeed, making circles $a$ or $b$ slightly bigger would turn the dotted line into a disk,
increasing the value of the depth of the tricycle $(1,a,b)$.
Using \eqref{eq:Descartes}, we can write this condition as 
$$
(a+b+1)^2 \ \geq \ 2(a^2 + b^2 + 1^2)\,.
$$
With a little algebra, this can be rewritten (for the equal sign) as
\begin{equation}
\label{eq:parabola}
2(a+b) = (a-b)^2 +1\,,
\end{equation}
which is evidently a parabolic equation.
Replace $x=a-b$, $y=a+b$ to get $2y=x^2+1$.
The main axis of the parabola coincides with the diagonal of the chart.
\\

Here is a recreational problem:
\\

{\small
\noindent
{\bf Problem:}  What is the probability that three random circles can be encompassed 
by a circle when put in a mutually tangent position?
\\[5pt]
{\bf Solution:}  
Clearly, the answer depends on the choice of measure one imposes on the set of circles.
Let us chose the uniform measure on the parameter of curvature.  
Using the above result, via a simple integral one finds the area of the parabolic region  
to make 3/4 of the square.
Thus the probability is $P = .75$.
}

%
%
%
%

\paragraph{B.  The elliptic plateau  of depth 2.}
We shall show that this region has an elliptic shape.
Figure \ref{fig:shape}(right) shows
a configuration in which disks $(1,a,b)$ 
are separated from the outer disk of curvature 0 (dotted line) by one 
intermediate disk of curvature $s$,
hence it is of depth 2.
The depth value is unstable: were the disk $s$ slightly smaller, 
the dotted line would become an inner disk, and   
the depth of the configuration would increase.
Thus the boundary of the region of depth 2 corresponds to 
this type of configuration.
The implied equations are (cf., \eqref{eq:Descartespm}): 
$$
\begin{array}{l}
s\ = \ a+b+0+2\sqrt{ab}\\
s\ = \ a+b+1-2\sqrt{ab+a+b}\\
\end{array}
$$
Eliminating $s$, we get 
$$
\sqrt{ab}+\sqrt{ab+a+b}\ = \ \frac{1}{2}\,,
$$
which after being squared twice, leads to
$$
(a+b)^2 -\tfrac{1}{2}(a+b) -ab +\frac{1}{16} \ = \ 0
$$
Substitution $x=a-b$ and  $y=a+b$ gives the standard form:
$$
                        12 x^2 + 36\left(y-\tfrac{1}{3}\right)^2 = 1
$$
which describes an ellipse.  In terms of the original variables:
\begin{equation}
\label{eq:ellipse}
36\left(a+b-\tfrac{1}{3}\right)^2 + 12 (a-b)^2 = 1
\end{equation}

Finding equations of the consecutive regions becomes increasingly more complex.
The next example shows more clearly the general method
of obtaining the equations for for other ellipses.

\begin{figure}
\centering
\includegraphics[scale=.77]{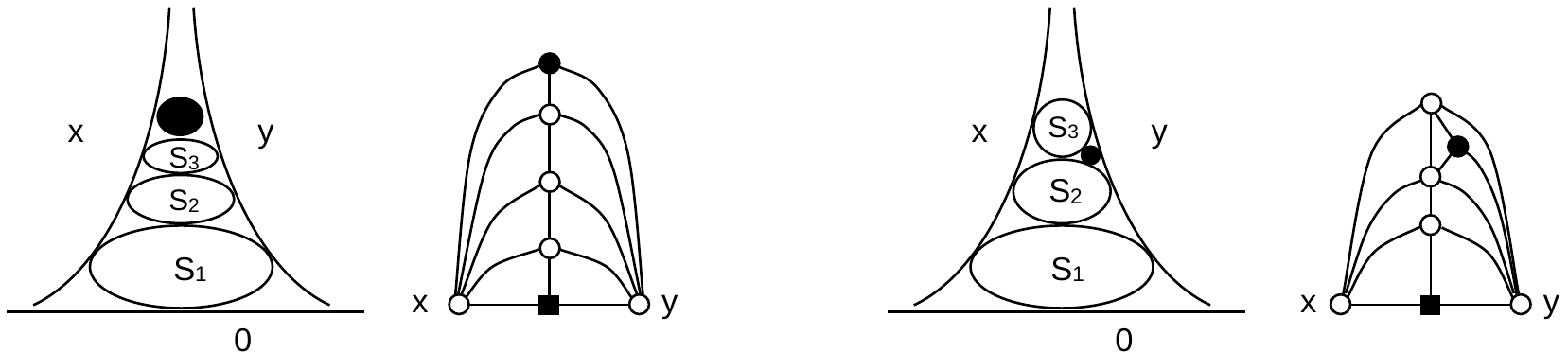}
\caption{The templates for two regions $\delta=4$.  
Disks in the arrangements are deformed in order to improve visualization.
Each arrangement is accompanied by a graph: black vertex denotes the disk of curvature 1, 
the square vertex stands for the straight line.}
\label{fig:d4}
\end{figure}

\paragraph{C.  The plateaus  of depth 4.}
Starting with he previous disk arrangements, we may increase the number of the disks 
separating the disk 1 from the ``tipping' dotted line.
We shall see two examples of such arrangements for regions of $\delta=4$,
see Figure \ref{fig:d4}. 

~

\noindent
{\bf Case 1:}
The idea is to write one quadratic equation, namely the Descartes formula, 
for a chosen Descartes configuration in this system, and write a 
set of linearized version \eqref{eq:Descarteslinear} for the remaining chain of disks.
Let us choose the bottom disks that includes the line (disk of curvature 0)
for the quadratic equation.
We get a system:
$$
\begin{array}{rll}
(A) &  (0+s_1 + x+y)^2 = 2(0^2+s_1^2+x^2 + y^2) \\[3pt]
(B) & 
\left\{
\begin{array}{ccc}
0+s_2&=&2(x+y+s_1)\\
s_1+s_3&=&2(x+y+s_2)\\
s_2+1&=&2(x+y+s_3)
\end{array}
\right.\\
\end{array}
$$
Part (B) consists of 3 linear equations, sufficient to express every $s_i$ in terms of $x$ and $y$.
Here we need $s_1=1/4 - 3x-3y$.
Under substitution, (A) becomes a quadratic equation in  $x$ and $y$:
$$
16x^2 +28xy+16y^2 -2x-2y +\dfrac{1}{16}  \ = \  0  \,,
$$
which describes one of the ellipses bounding region od $\delta=4$,
namely the one in the diagonal chain (labeled later as $G4_{1,1}$). 

~\\
{\bf Case 1.5:}  The quadratic equation for the system of equations 
may be chosen also from a different Descartes configuration in the arrangement.
For instance in the above example we may consider the configuration involving disks $x$, $y$, and 1:
$$
(A) \quad 
(1+s_3 + x+y)^2 = 2(1^2+s_3^2+x^2 + y^2) 
$$
The 3 linear equations (B)  stay the same, but now we need to extract
a different unknown:  $s_3=3/4 - 3x-3y$.
Under substitution, the quadratic equation becomes:
$$
16x^2+28xy + 16y^2 -2x-2y+\frac{1}{16}=0
$$
as before.

~\\
{\bf Case 2:} Similarly, for the 
third arrangement in Figure \ref{fig:d4}, we have: 
$$
\begin{array}{rll}
(A) &
(1+ x+y+ s_2)^2 = 2(1^2+x^2 + y^2+s_2^2) \\[3pt]
(B) & 
\left\{
\begin{array}{ccc}
0+s_2&=&2(s_1+s_3+y)\\
s_1+x&=&2(s_2+s_3+y)\\
s_3+1&=&2(x+y+s_2)
\end{array}
\right.
\end{array}
$$
Substituting $s_2=2/5 - 2x/3 -6y/5$ extracted from (B) to equation (A), we get:
$$
\frac{25}{9}x^2+\frac{22}{3}xy + \frac{121}{25}y^2 -2x-\frac{34}{25}y+\frac{9}{25}=0
$$

\paragraph{D.  The general method.}
To obtain other equations,
one needs to consider the variety of possible disk pyramidal arrangements,
by which we mean the following:
Set a line (disk of zero curvature) at the bottom.
Set two disks mutually tangent, call it A and B.
In the ideal triangle formed  in the space 
enclosed, inscribe a system of disks
so that ...
The smallest is of curvature $1$.
Choose two of the disks tangent to $1$ and  denote curvatures $x$ and $y$.
Denote the remaining curvatures by $s_1$, $s_2$,...,$s_n$.
Now write a system of equations:
(A) a quadratic equation for one of the Descartes configurations, 
and a system of linear equations (B)
for the remaining chain of disks using 
\eqref{eq:Descarteslinear}.

%
%

\paragraph{E. Labeling system}

There is a problem of uniquely labeling the ellipses in the spiderweb.
One way is to follow the ternary tree structure. 
E.g., counting from say 2,  2DRRL means: move from 2 first down one step, 
two steps to the right and one step to the left, every time decreasing the depth by 1, where L, R, D
need to be defined e.g, as ``left'', ``right'', and ``down''
with respect to the direction from the last entry.

~

\begin{figure}[H]
\includegraphics[scale=.83]{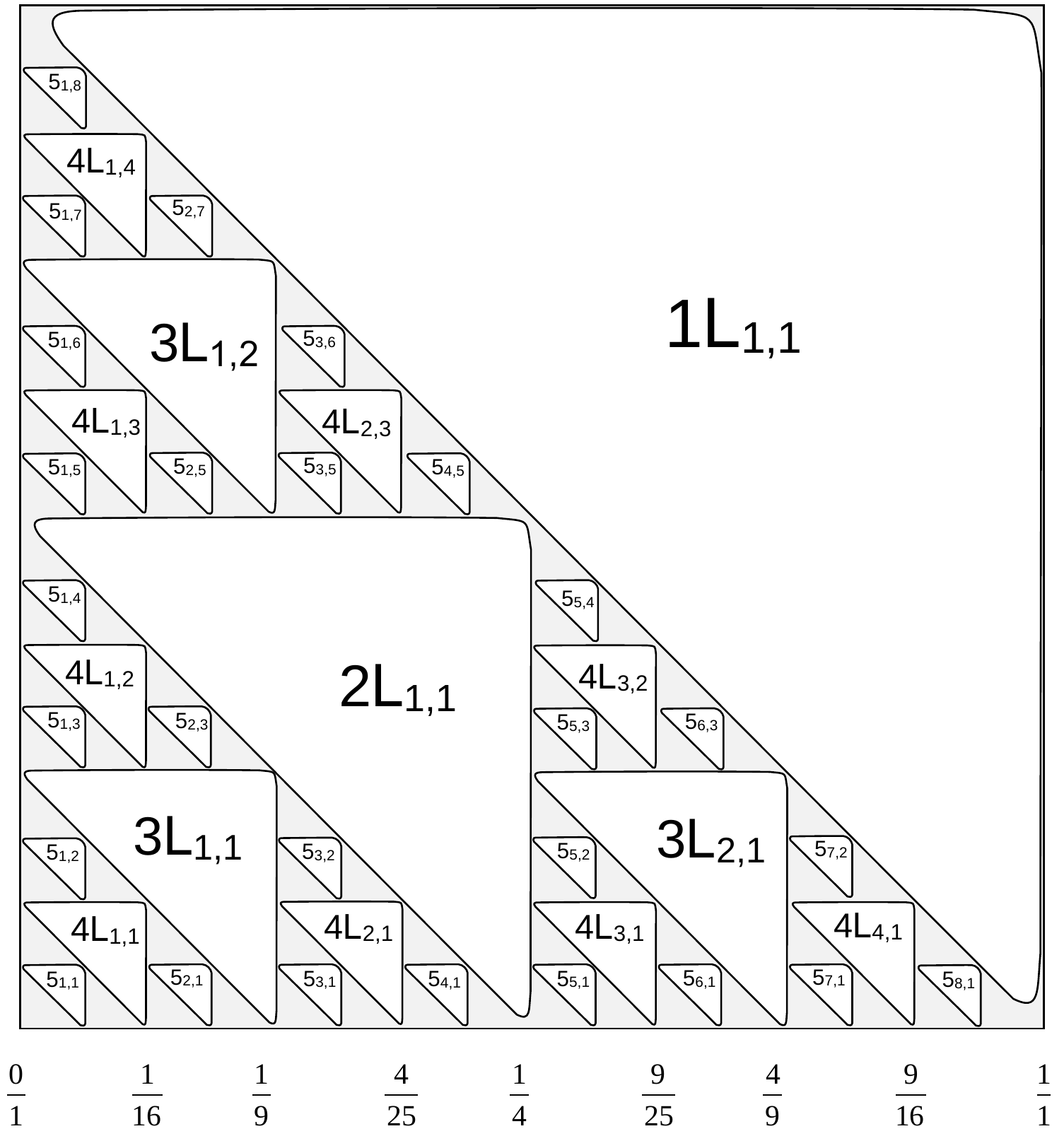}
\caption{A template for labeling of the regions.
Letter ``L'' is suppressed in smaller triangles, thus $5_{3,1}\equiv 5L_{3,1}$.
In this graphic representation, the size of a triangle corresponds to the value of the depth function.}
\label{fig:Sierpinski1}
\end{figure}

Instead, we use the following labeling system:
Deform the web so that the ellipses become triangles 
organized in an orderly fashion, as shown in Figure \ref{fig:Sierpinski1}.
Note the benefits:
\begin{enumerate}
\item
Regions of the same value of the Apollonian depth function $\delta$ become congruent triangles.
\item
For a fixed value of $\delta$ the corresponding triangles form a matrix-like pattern.
\end{enumerate}
This allows one to use labeling:
$$
dL_{i,j}
$$
where ``L'' is just a standard indication that this particular labeling is used,
``$d$'' is the value of the depth,
and $i,j$ are the ``coordinates'' of the triangle within this set of triangles:  
$i$ stands for the number of the column
(from the left).
and $j$ stands for the row (starting at the bottom).
Figure \ref{fig:Sierpinski1} shows the labeling with space-saving omission of the letter ``L''.

Note that topologically the pattern is that of Sierpi\'nski triangle.


\paragraph{F. Quadratic equations.}

Here are the equations for the first few plateaus, obtained by the method outlined above. 
$$
\begin{array}{lll}
1L_{1,1}:\quad
& x^2-2xy+y^2-2x-2y+1    \ = \  0   \\[10pt]
2L_{1,1} :
& 4x^2+4xy+4y^2-2x-2y+1/4    \ = \  0  \\[10pt]
3L_{1,1}:
& 9x^2+14xy+9y^2 -2x-2y+ \dfrac{1}{9}  \ = \  0  \\[10pt]
3L_{2,1}:
& \dfrac{9}{4} x^2+2xy +4y^2-2x-  \dfrac{4}{3} y  +\dfrac{4}{9}  \ = \  0  \\[10pt]
4L_{1,1}:
& 16x^2 +28xy+16y^2 -2x-2y +\dfrac{1}{16}  \ = \  0  \\[10pt]
4L_{3,2}:
&  \dfrac{196}{81}x^2+  \dfrac{46}{27} xy +\dfrac{121}{36} y^2 
                       - \dfrac{52}{27} x - \dfrac{14}{9} y+ \dfrac{4}{9}    \ = \  0  \\[10pt]
\end{array}
$$
A few curious things become at once visible.
\begin{enumerate}
\item
Often, but not always, the free term is reciprocal to the coefficient at $x^2$.
In such a case, the free term is the x-coordinate of the point the ellipse touches the $x$-axis.
\item
 If the ellipse is tangent to the $x$-axis, its coefficient at the term $x$ is equal to 2.
 \item
Another intriguing property is the following: completing the square
leaves out the same term, $2xy$ in each of the equations.

\end{enumerate}

The above equations may be thus rewritten with the squares completed
and in a form that should make inspection easier: 
$$
\begin{array}{lllll}
1L_{1,1}: \quad &  \left( x+y \right)^2 +1    &\quad = \quad &4xy + 2x+2y \\[10pt]
2L_{1,1}:  &  \left( 2x+2y \right)^2+\dfrac{1}{4}    &\quad = \quad  &4xy + 2x+2y \\[10pt]
3L_{1,1}:  &   \left( 3x+3y\right)^2  + \dfrac{1}{9}  &\quad = \quad  &4xy + 2x+2y \\[10pt]
3L_{2,1}:  &  \left( \dfrac{3}{2} x +2y \right)^2   +\dfrac{4}{9}  &\quad = \quad  &4xy + 2x +  \dfrac{4}{3} y \\[10pt]
4L_{4,1}:   &  \left( 4x +4y \right)^2 +\dfrac{1}{16}  &\quad = \quad  &4xy + 2x+2y  \\[10pt]
1L_{3,2}:   &  \left(  \dfrac{14}{9}x + \dfrac{11}{6} y \right)^2  + \dfrac{4}{9}    
                                &\quad = \quad  &4xy +  \dfrac{52}{27} x + \dfrac{14}{9} y  \\[10pt]
%
%
\end{array}
$$

A longer list of the ellipse equations may be found in Appendix B.

The above properties suggested the theorem that is spelled out in the next section. 

\newpage

\section{Stern-Brocot structure}

The most  intriguing property of the fractal is appearance 
of Stern-Brocot stucture in the pattern of the points of tangency.

\subsection{The x-axis corona} 
  
Inspect the $x$-wing main chain, the chain of ellipses of depth 2,3,4, etc., that  
are simultaneously tangent to the parabolic region and the $x$-axis.
Using the equations of the previous section
we may find that the $x$-coordinates of the points 
they touch the axis
form  the following progression of fractions:
$$
\frac{1}{4}, \quad
\frac{4}{9}, \quad
\frac{9}{16}, \quad 
\frac{16}{25}, \quad ..., \qquad
\frac{(n-1)^2}{n^2}
$$

\begin{figure}[h]
\centering
\includegraphics[scale=.78]{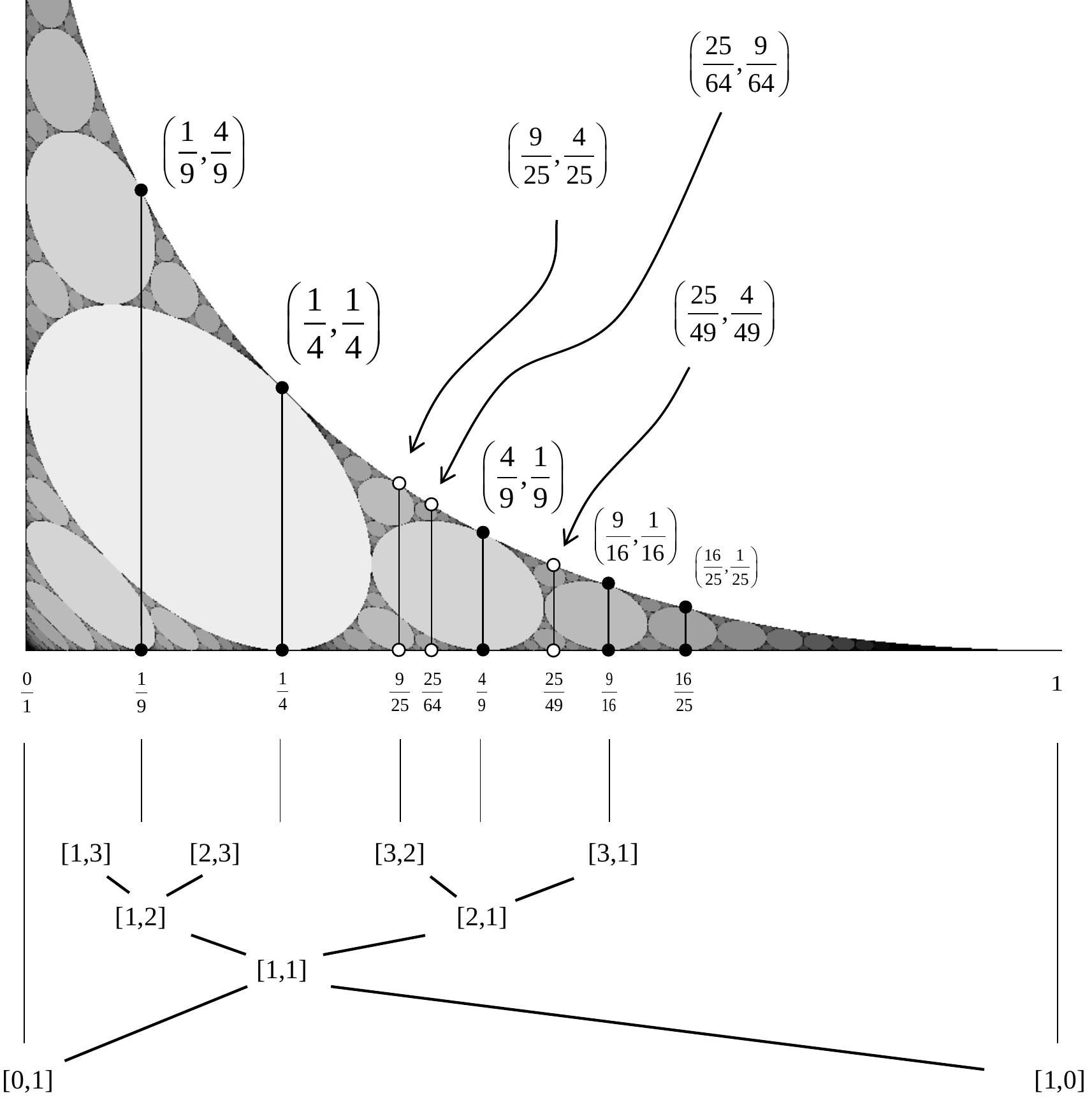} 
\caption{\small The spiderweb fractal and the Stern-Brocot tree}
\label{fig:chart-tree2}
\end{figure}

\begin{figure}
\includegraphics[scale=.8]{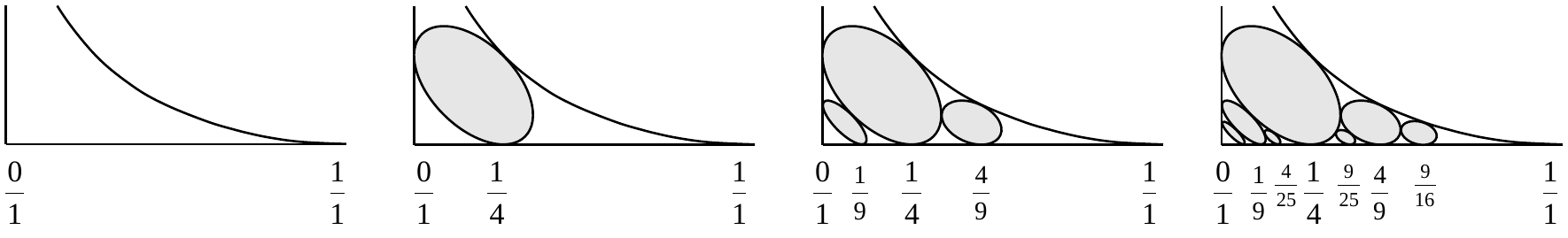}
\caption{Producing the x-points}
\label{fig:x-comix}
\end{figure}

Further inspection reveals that 
the ellipses inscribed between them touch the $x$-axis at points that 
result via a 
deformed version of the Farey addition of fractions.
Namely, 
the point for the inscribed ellipse between
$\frac{a}{b}$ and $\frac{c}{d}$
is 
\begin{equation}
\label{eq:myFarey}
\frac{a}{b}\boxplus\frac{c}{d} 
\ = \ \frac{(\sqrt{a}\; +\sqrt{c})^2}{(\sqrt{b}\; +\sqrt{d})^2}
\end{equation}
Iterating this process will account for all ellipses tangent to the $x$-axis, 
i.e., the $x$-corona.
Figure \ref{fig:x-comix} shows an order of recovering the points,
which starts with the extreme fractions 0/1 and 1/1,
and then follows the deformed Farey addition.

~

Before we collect these facts in one extended statement, let us recall the 
basic facts of a  Stern-Brocot tree.
\\

%

\subsection{Definition of the Stern-Brocot array} 

Start with a 2-element sequence of pairs:  [1,0], [0,1] (vectors).
Create an array of sequences, making new sequences from the previous by 
inscribing new terms between the existing terms.
The new terms are simply the vector sums of the neighbors:
\begin{equation}
\label{eq:abcd}
          ...,\ [a,b],\ [c,d],\ ... \qquad \mapsto \qquad   ...,\ [a,b],\ [a\!+\!c,b\!+\!d], \ [c,d],\ ... 
\end{equation}
The result should be called the {\bf Stern-Brocot array}, the initial fragment is presented below:
$$
\begin{array}{ccccccccccc}
[1,0]   &&&         &&         &&&    [0,1] \\[7pt]
[1,0]   &&&&       [1,1]    &&&&    [0,1] \\[7pt]
[1,0]   &&[2,1]&& [1,1]    &&[1,2]&&     [0,1] \\[7pt]
[1,0] \   &\ [3,1] \ &\ [2,1] \ & \ [3,2] \ &\  [1,1] \   & \ [2,3]\ &\ [1,2]\ &\ [1,3]\ &  \   [0,1] \\
\end{array}
$$

\newpage
\noindent
By removing the multiple occurrence of terms, the array becomes a{\bf tree}:
\begin{center}\label{stern-brocot}
\begin{tikzpicture}[level distance=0.5in]
\tikzstyle{level 1}=[sibling distance=2in]
\tikzstyle{level 2}=[sibling distance=1in]
\tikzstyle{level 3}=[sibling distance=0.5in]
\node {$[1,1]$}
  child {node {$[1,2]$}
    child {node {$[1,3]$}
      child {node {$[1,4]$}
      }
      child {node {$[4,3]$}
      }
    }
    child {node {$[3,2]$}
      child {node {$[3,5]$}
      }
      child {node {$[5,2]$}
      }
    }
  }
  child {node {$[2,1]$}
    child {node {$[2,3]$}
      child {node {$[2,5]$}
      }
      child {node {$[5,2]$}
      }
    }
    child {node {$[3,1]$}
      child {node {$[3,4]$}
      }
      child {node {$[4,1]$}
      }
    }
  };
\end{tikzpicture}
\end{center}
The main property of the Stern-Brocot tree is that the pairs are relative primes.
In particular, by replacing
$$
[p,q] \quad \mapsto \quad \frac{p}{q}\,,
$$
the tree becomes a tree of all positive rational numbers,
and this is how the original tree is usually defined and presented.
In such a case, the the rule \eqref{eq:abcd} is replaced by the so-called
Farey  addition of fractions:
$$ 
 \frac{p}{q}\oplus \frac{s}{t} \ = \ \frac{p+s}{q+t}
$$

\newpage
\subsection{The main statement on ellipses and squared fractions}

%
%
%
%

%
%
%
%
%
%

The following property is analogous to the well-known Ford's theorem for circles \cite{Ford}.

~

\begin{figure}[h]
\centering
\includegraphics[scale=.7]{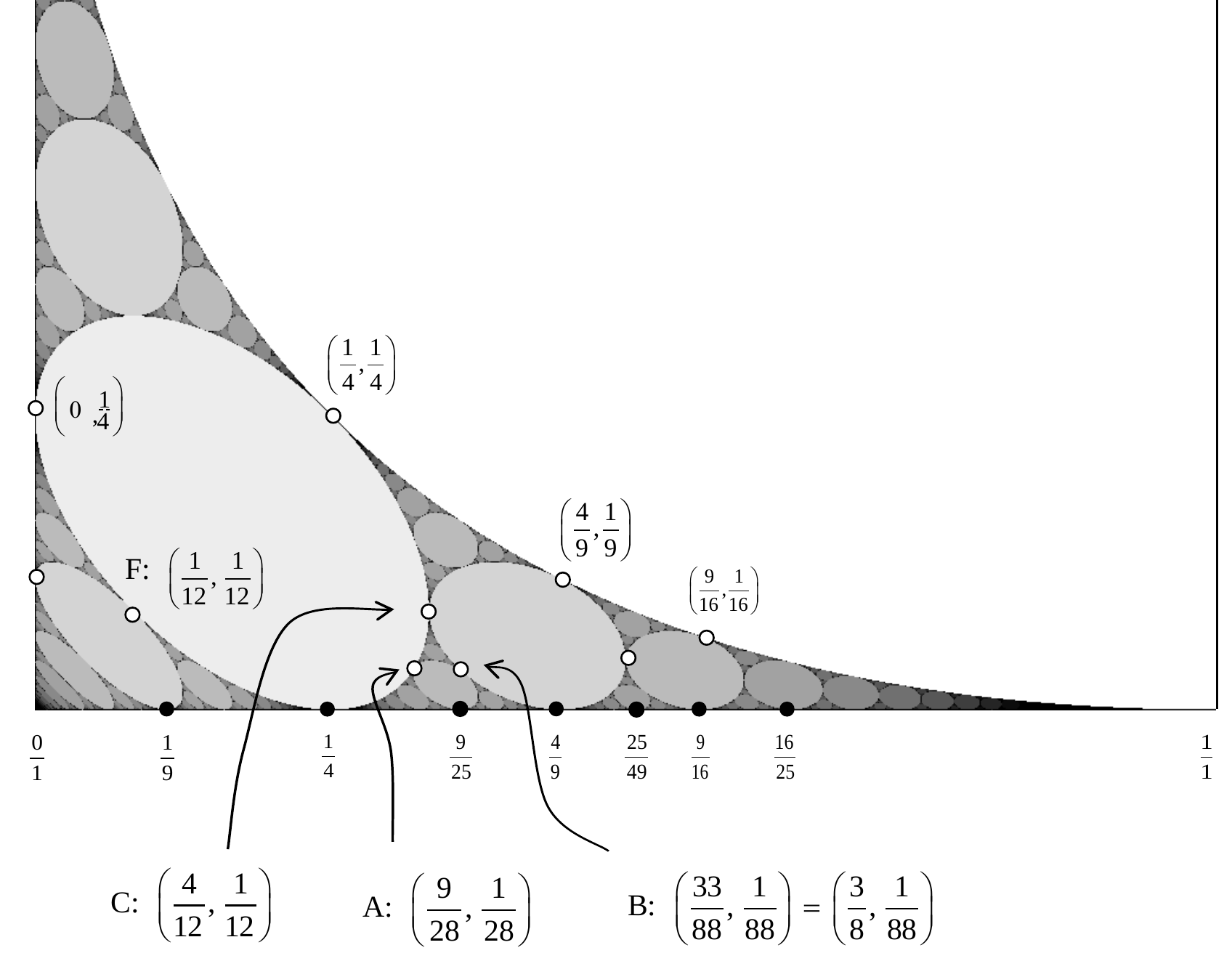}
\caption{\small Illustration of Theorem 1(C) }
\label{fig:Theorem}
\end{figure}

\noindent
{\bf Theorem 1:}
{\sf 
For every reduced rational number $\frac{p^2}{m^2}$ 
in the closed interval $[0,1] \subset \mathbb R$ draw an ellipse $E[p,m]$ defined by  
\begin{equation}
\label{eq:my}
\left(\,\frac{m}{p}\;\mathbi x + \frac{p^2\!+\! m^2\! -\! 1}{pm} \; \mathbi y\,\right)^2 + \frac{p^2}{m^2}
 \ = \ 
4\mathbi x \mathbi y + 2\mathbi x + 2\,\frac{m^2\!-\!p^2\! +\! 1}{m^2}\,\mathbi y
\end{equation}
Then the following hold:
\\[7pt]
{\bf (A)} Each ellipse lies above the $x$-axis and is tangent to it at the point $p^2/m^2$. 
\\[7pt]
{\bf (B)} The interiors of the ellipses are disjoint.
Moreover, if fractions $\frac{p}{m}$ and $\frac{q}{n}$ 
satisfy
\begin{equation}
\label{eq:det}
\det \begin{bmatrix}p &q \\ m& n\end{bmatrix} \ =\ \pm1
\end{equation}
then there is an ellipse inscribed between these ellipses and the $x$-axis, namely
 ellipse $E[p\!+\!q,m\!+\!n]$ over the point 
\begin{equation}
\label{eq:x}
x=\frac{(p\!+\!q)^2}{(m\!+\!n)^2}\,.
\end{equation} 
{\bf (C)} Two ellipses satisfying \eqref{eq:det} are mutually tangent at point:
\begin{equation}
\label{eq:boki}
E[p,m] \cap E[q,n] \ \ = \ \ \ \left(\,\frac{p^2+q^2-1}{m^2+n^2-1},  \; \frac{1}{m^2+n^2-1} \,\right) \,.
\end{equation}
~\\[7pt]
{\bf (D)} The ellipses coincide with the plateaus of the constant values of the depth function $\delta$.
In particular,
the value of $\delta$ in $E[p,n]$ is equal to the row number of $(p,n)$ in the 
Stern-Brocot tree X.
}

~\\
{\bf Proof:}  Lengthy  calculations with a support of computer simulations..
\qed
\\

\begin{figure}[H]
\centering
\includegraphics[scale=.68]{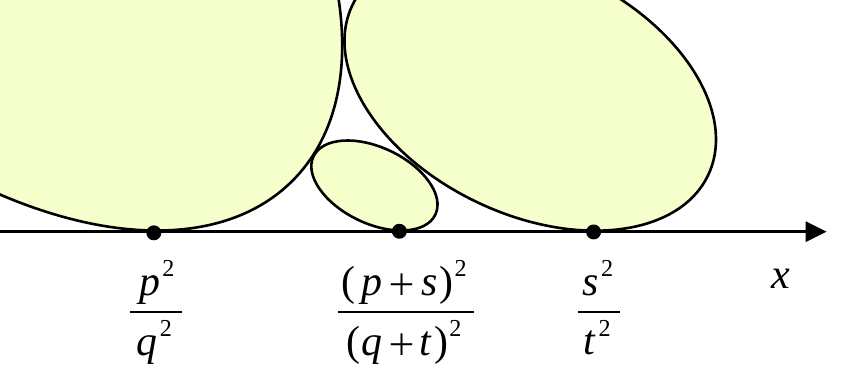}
\quad\includegraphics[scale=.68]{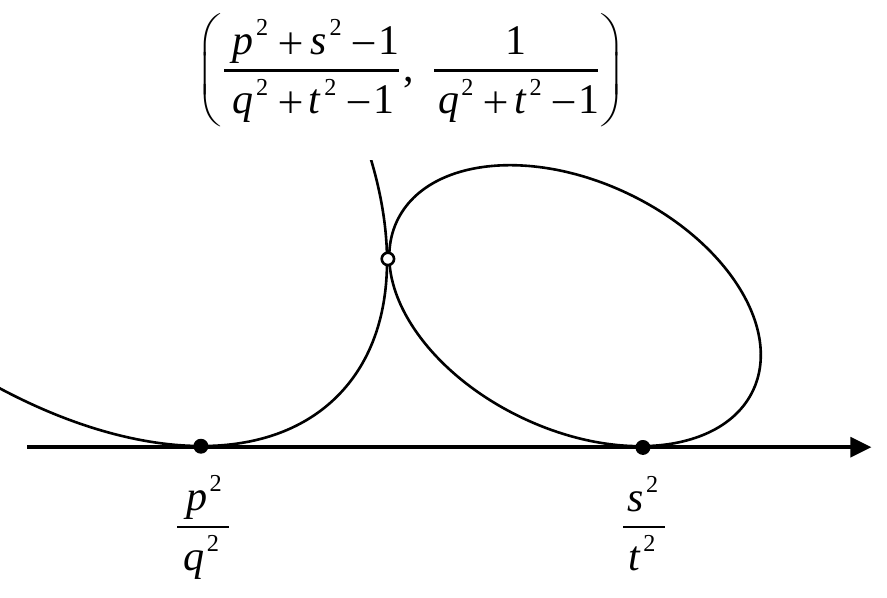}
\caption{\small Mid-points.  Left:  at the bottom (super-axial);  Right: sub-parabolic}
\label{fig:}
\end{figure}

\noindent
{\bf Special case 1:}
Note: the main chain of ellipses is formed above the fractions such that $m=p+1$.
The equations of the ellipses have a detectable pattern
$$
\left(\,\frac{p+1}{p}\;\mathbi x + 2\; \mathbi y\,\right)^2 + \frac{p^2}{(p+1)^2} \ = \ 
4\mathbi x \mathbi y + 2\mathbi x + \frac{4}{p+1}\,\mathbi y
$$

\noindent
{\bf Special case 2:}
The parabola is included in the $x$-corona, namely as $E[1,1]$
tangent to $x$-axis at $x=(1/1)^2=1$.  
Indeed, the general equation \eqref{eq:my} reduces to the parabola equation $2(x+y)=(x-y)^2+1$
under substitution  $p=1$ and $m=1$. 
The points of tangency along the parabola with the main chain are now the points described by \eqref{eq:boki}.
Consult Figure \ref{fig:Theorem}.

~

\noindent
{\bf Special case 3:}
The vertical line at $x=0$ should also be considered as a part of the x-corona, namely as $E[0,1]$.,
Substitution $p=0$ and $m=1$ reduces \eqref{eq:my} to $x=0$ 
(to avoid infinities, one must multiply both sides of X by $p^2$ before the substitutions). 
The theorem holds with the inclusion of this line. 
For instance,  the points of the tangency of the the ellipses of the diagonal chain with the $y$-axis
agree with Equation \eqref{eq:boki} and are:
$$
E[1,m] \cap E[0,1] \ \ = \ \ \ \left(\,\frac{0}{1},  \; \frac{1}{m^2} \,\right) \,:\quad
\left(\frac{0}{1},\, \frac{1}{4} \right), \quad
\left(\frac{0}{1},\, \frac{1}{9} \right), \quad
\left(\frac{0}{1},\, \frac{1}{16} \right), \quad
etc.
$$

~

\noindent
{\bf Special case 4:}
Each ellipse in the main $x$-wing chain is tangent to the $x$-axis
and to the parabola at points that are vertically aligned.
The coordinates on the parabola are 
$$
E[p,p\!+\!1] \cap E[1,1] \ \ = \ \ \ \left(\,\frac{p^2}{(p\!+\!1)^2},  \; \frac{1}{(p\!+\!1)^2} \,\right) \,.
$$

~

The statement \eqref{eq:x}
implies the following structure of the tangency points for the $x$-corona:
\begin{center}\label{stern-brocot}
\begin{tikzpicture}[level distance=0.5in]
\tikzstyle{level 1}=[sibling distance=2.4in]
\tikzstyle{level 2}=[sibling distance=1.2in]
\tikzstyle{level 3}=[sibling distance=0.5in]
\node {$\frac{1}{4}$}
  child {node {$\frac{1}{9}$}
    child {node {$\frac{1}{16}$}
      child {node {$\frac{1}{25}$}
      }
      child {node {$\frac{4}{49}$}
      }
    }
    child {node {$\frac{4}{25}$}
      child {node {$\frac{9}{64}$}
      }
      child {node {$\frac{9}{49}$}
      }
    }
  }
  child {node {$\frac{4}{9}$}
    child {node {$\frac{9}{25}$}
      child {node {$\frac{16}{49}$}
      }
      child {node {$\frac{25}{64}$}
      }
    }
    child {node {$\frac{9}{16}$}
      child {node {$\frac{25}{49}$}
      }
      child {node {$\frac{16}{25}$}
      }
    }
  };
\end{tikzpicture}
\end{center}

~

\subsection{Parabolic corona}

The parabolic corona
consists of the ellipses tangent to the parabolic region of depth $\delta =1$.
It includes the x-wing main chain, the $y$-wing main chain,
and the ellipses inscribed in the regions between them and the parabola.
It turns out that 
the pattern of the tangency follows the same Stern-Brocot structure with the deformed Farey addition, 
but now it applies to both coordinates $x$ and $y$, and extends over the whole length of the parabola.

~

\noindent
{\bf Proposition}
{\sf  For any positive value of $p$ and $q$ and $n=p+q$,
the numbers 
$$
a = \frac{p^2}{n^2},\quad b= \frac{q^2}{n^2}
$$ 
satisfy the parabolic Equation \eqref{eq:parabola},  
thus the points $(a,b)$ lie on the parabola $P$.
}

~

In this context, the pairs $[p,q]$ will 
label both the points on the parabola with coordinates 
 $$
         \mathbb Z^2  \ni [p,q] \quad \mapsto\quad (x,y)\ =\ \left(\frac{p^2}{n^2},\; \frac{q^2}{n^2}\right)  \in \mathbb R^2
      \qquad \hbox{where}\quad    n=p+q
$$      
and ellipses tangent to the parabola at that points.            
Such ellipses will be denoted by $F[p,q]$.

~

The parabolic corona undergoes a phenomenon 
analogous to that of the $x$-axis corona:

~

\noindent
{\bf Proposition:}  
{\sf The ellipses in the parabolic corona in the web 
are tangent to the parabola in the 
rational points of form
$$
(x,y) \ = \left(   \frac{p^2}{n^2},\; \frac{q^2}{n^2}\right), \qquad n=p+q\,.
$$ 
{\bf (A)} For every $(p,q)\in\mathbf N$ there is such an ellipse $F[p,q]$ 
tangent to the parabola at the above point.
Additionally, two ellipses, $F[p,q]$ and $F[p',q']$
are mutually tangent iff   
$$
\det\begin{bmatrix} p&p'\\ m&m'\end{bmatrix}=\pm 1  
\quad\hbox{and}\quad
\det\begin{bmatrix} q&q'\\n&n'\end{bmatrix}=\pm 1  
\qquad m=p\!+\!q,\ m'=p'\!+\!q'
\,.
$$
(The two conditions are equivalent due to \eqref{eq:parabola}.) \; The point of tangency is 
$$
\left(   \frac{p^2+p'^2-1}{m^2+n^2-1},\ \frac{q^2+q'^2-1}{m^2+n^2-1} \right)
$$
{\bf (B)} The ellipse inscribed between the above ellipses is
$$
F[p+p',q+q']\,,
$$
touching the parabola at the point with coordinates resulting from the deformed 
Farey addition:
$$
\left(\frac{p^2}{m^2},\,\frac{q^2}{m^2}\right)  \boxplus  \left(\frac{p'^2}{m'^2},\,\frac{q'^2}{m'^2}\right)
\ = \ 
\left(\frac{(p\!+\!p')^2}{(m\!+\!m')^2} ,\; \frac{(q\!+\!q')^2}{(m\!+\!m')^2} \right)
$$

}

~

~

\begin{figure}[h]
\centering
\includegraphics[scale=.68]{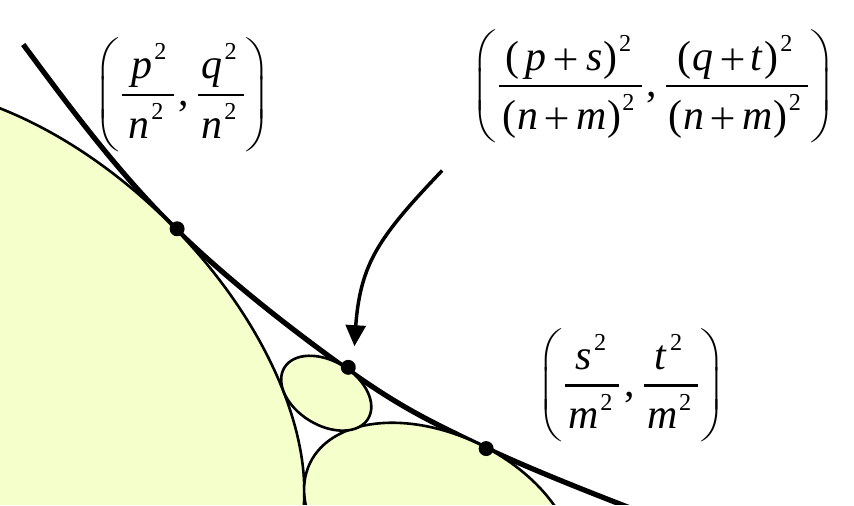}
\qquad
\includegraphics[scale=.68]{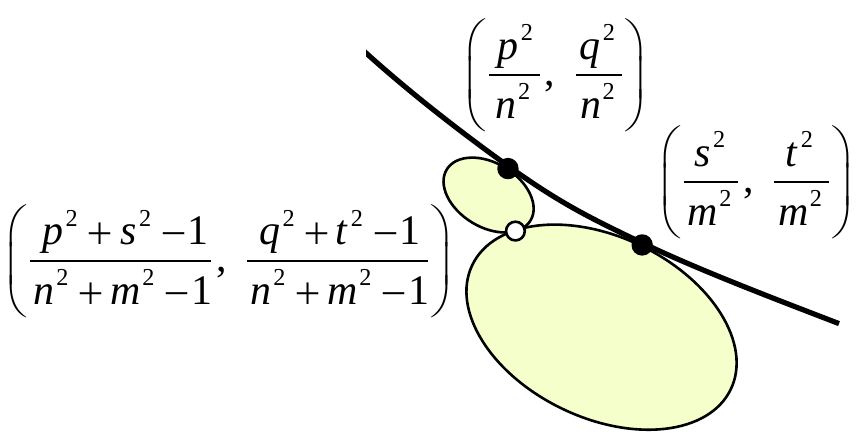}
\caption{\small Tangential points.  Left:  between ellipses and the parabola; 
 Right: between tangent ellipses}
\label{fig:thm2}
\end{figure}

Figure \ref{fig:thm2} summarizes the main points.

~

Note that as before, one may organize all ellipses in the parabolic corona and the corresponding points on the parabola
in a form of a tree.

\begin{center}\label{stern-brocot}
\begin{tikzpicture}[level distance=0.5in]
\tikzstyle{level 1}=[sibling distance=2.7in]
\tikzstyle{level 2}=[sibling distance=1.4in]
\tikzstyle{level 3}=[sibling distance=0.7in]
\node {$\left(\frac{1}{4},\,\frac{1}{4}\right)$}
  child {node {$\left(\frac{1}{9},\,\frac{4}{9}\right)$}
    child {node {$\left(\frac{1}{16},\,\frac{9}{16}\right)$}
      child {node {$\left(\frac{1}{25},\,\frac{16}{25}\right)$}
      }
      child {node {$\left(\frac{4}{49},\,\frac{25}{49}\right)$}
      }
    }
    child {node {$\left(\frac{4}{25},\,\frac{9}{25}\right)$}
      child {node {$\left(\frac{9}{64},\,\frac{25}{64}\right)$}
      }
      child {node {$\left(\frac{9}{49},\,\frac{16}{49}\right)$}
      }
    }
  }
  child {node {$\left(\frac{4}{9},\,\frac{1}{9}\right)$}
    child {node {$\left(\frac{9}{25},\,\frac{4}{25}\right)$}
      child {node {$\left(\frac{16}{49},\,\frac{9}{49}\right)$}
      }
      child {node {$\left(\frac{25}{64},\,\frac{9}{64}\right)$}
      }
    }
    child {node {$\left(\frac{9}{16},\,\frac{1}{16}\right)$}
      child {node {$\left(\frac{25}{49},\,\frac{4}{49}\right)$}
      }
      child {node {$\left(\frac{16}{25},\,\frac{1}{25}\right)$}
      }
    }
  };
\end{tikzpicture}
\end{center}

For instance
$$
\left(\frac{1}{4},\frac{1}{4}\right)  \boxplus  \left(\frac{4}{9},\frac{1}{9}\right)
=\left(\frac{1}{4}\boxplus\frac{4}{9} ,\; \frac{1}{4}  \boxplus \frac{4}{9}\right)
=\left(\frac{(1+2)^2}{(2+3)^2} ,\; \frac{(1+1)^2}{(2+3)^2} \right)
=\left(\frac{9}{25} ,\; \frac{4}{25} \right)
$$

\begin{figure}[H]
\includegraphics[scale=.77]{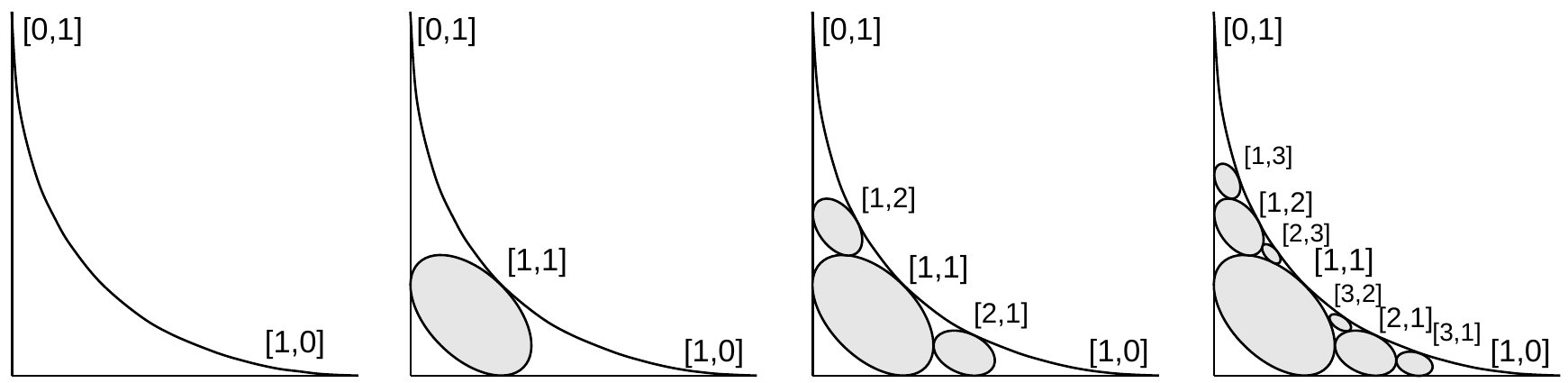}
\caption{Producing the P-points}
\label{fig:p-comix}
\end{figure}

Figure \ref{fig:p-comix}
shows the first few steps of such recurrence.
Note that the y-axis and the x-axis are among the ellipses
as the special cases $F[0,1$ and $F[1,0]$. 

~

\noindent
{\bf Pseudoproofs:}
The claims were tested by a computer process of  evaluating
the depth function at the neighborhoods of the tangency points by probing such as    
$$
\delta (a,b) \ \hbox{and} \  \delta(a,b+\varepsilon) 
$$
With Maple,  the value of $\varepsilon = 0.0000001$ already leads to a jump.

%

~

\noindent
{\bf Remark:} 
The base points of the lower corona and the ceiling points of the upper corona are vertically aligned, 
i.e., they share the same $x$-coordinate. 

~

\newpage

\section{Barycentric coordinates}
\label{sec:barycentric}

\begin{figure}[h]
\centering
\includegraphics[scale=.45]{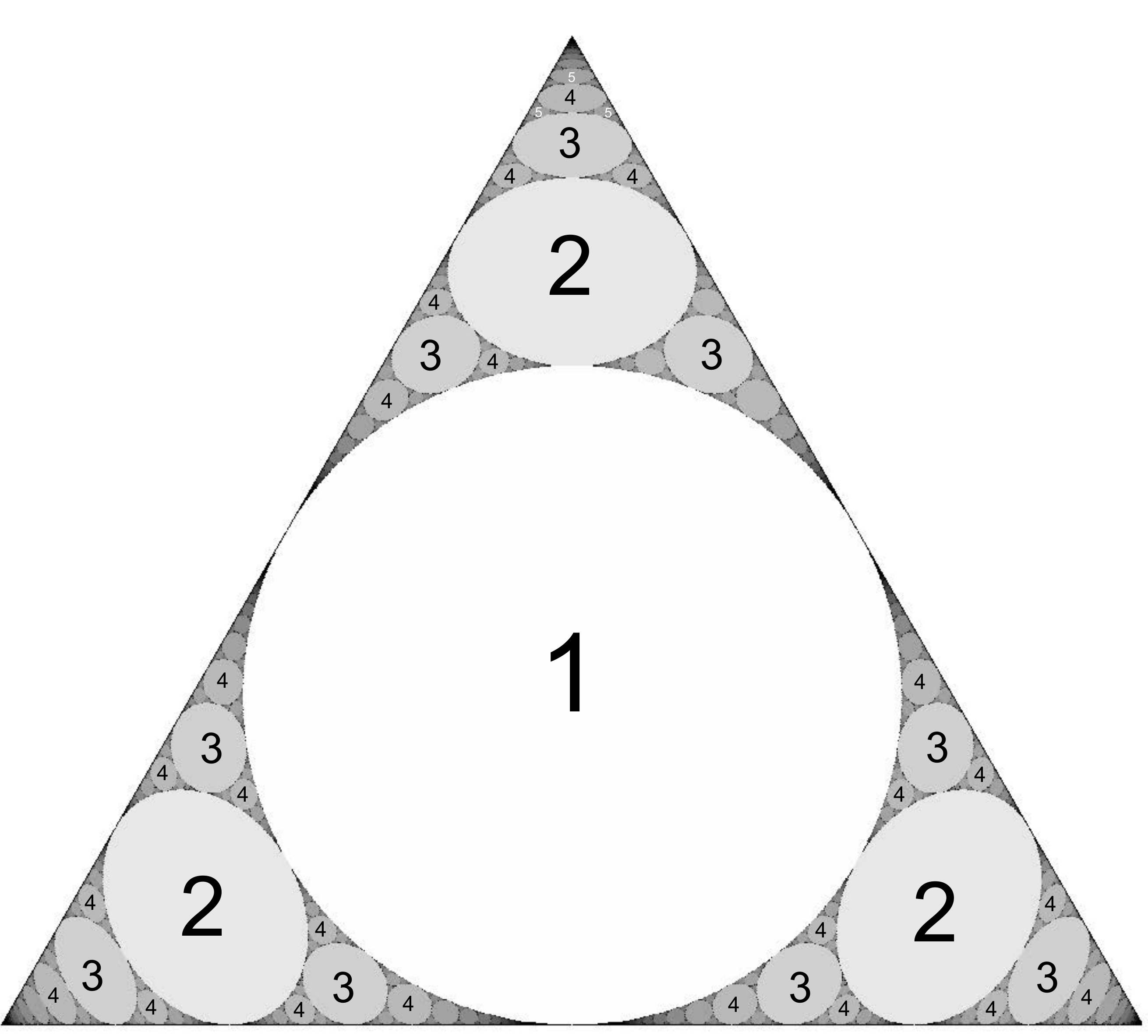} 
\caption{\small Depth in Barycentric representation. Figure obtained }
\label{fig:bary-nrs}
\end{figure}

As mentioned in Section \ref{sec:d},
 the chart of the depth function may be represented in barycentric coordinates.
We denote them with double brackets. 
The curvatures of a tricycle $(a,b,c)$ of disks are rescaled to 
$$
(\!(\,  x,\, y,\, z\,)\!) \ \ = \ \ \left( \frac{a}{a+b+c},\, \frac{b}{a+b+c},\, \frac{c}{a+b+c}\,  \right)
$$
($x+y+z=1$).
The resulting fractal is shown in Figure \ref{fig:bary-nrs}.

\noindent
In this setup the three disks of the tricycles are distinguishable.  
The three vertices of the triangle correspond to two lines separated by a circle, 
each time a different distribution:  $(1,0,0),\; (0,1,0),\; (0,0,1)$.
The center of the triangle corresponds to three equal in size disks, $(1/3,1/3,1/3)$.
The Apollonian Window is determined by (among others) point
of barycentric coordinates $\left(\left(\frac{2}{7},\, \frac{2}{7},\, \frac{3}{7}\,\right)\right)$. 
\\
\\
This form of the Apollonian spiderweb is more symmetric and elegant,
but for experimentation, the rectangular framework of the previous sections 
 is easier to explore.
\\

Topologically, the result is equivalent to disk packing of triangle. 
It has also the structure of Sierpi\'nski triangle.

\begin{figure}[H]
\centering
\begin{tikzpicture}[level distance=1cm,
  level 1/.style={sibling distance=2.4cm},
  level 2/.style={sibling distance=.78cm},
  level 3/.style={sibling distance=.2cm},
   level 4/.style={sibling distance=.1cm}]
     \node {1}
    child {node {2}
                  child {node {3}
                         child {node {\footnotesize 4}}
                         child {node {\footnotesize 4}}
                         child {node{\footnotesize 4}}
     }
                  child {node {3}
                         child {node {\footnotesize 4}}
                         child {node{\footnotesize 4}}
                         child{node{\footnotesize 4}} 
    }
                 child {node {3}
                         child {node {\footnotesize 4}}
                         child {node{\footnotesize 4}}
                         child{node{\footnotesize 4}} 
      } }
    child {node {2}
                  child {node {3}
                         child {node {\footnotesize 4}}
                         child {node{\footnotesize 4}}
                         child{node{\footnotesize 4}} 
    }
      child {node {3}
                          child {node {\footnotesize 4}}
                         child {node{\footnotesize 4}}
                         child{node{\footnotesize 4}} 
    }
                child {node {3}
                         child {node {\footnotesize 4}}
                         child {node{\footnotesize 4}}
                         child{node{\footnotesize 4}} 
      } }      
    child {node {2}
                  child {node {3}
                         child {node {\footnotesize 4}}
                         child {node{\footnotesize 4}}
                         child{node{\footnotesize 4}} 
    }
      child {node {3}
                          child {node {\footnotesize 4}}
                         child {node{\footnotesize 4}}
                         child{node{\footnotesize 4}} 
    }
                child {node {3}
                         child {node {\footnotesize 4}}
                         child {node{\footnotesize 4}}
                         child{node{\footnotesize 4}} 
      } }  
      ;
\end{tikzpicture} \quad
\quad
\begin{tikzpicture}[scale=1.2]
\draw [thick] (-1.732,-1)--(1.732,-1)--(0,2)--(-1.732,-1);  
\draw [thick] (0,0) circle (1); 
\draw [thick] (-1.154,-.667) circle (.333); 
\draw [thick] (1.154,-.667) circle (.333); 
\draw [thick] (0, 1.333) circle (.33); 
\draw [thick] (-1.54,-.888) circle (.111); 
\draw [thick] (1.54,-.888) circle (.111); 
\draw [thick] (0,1.777) circle (.111); 
\draw [thick] (-.733,-.867) circle (.132); 
\draw [thick] (.733,-.867) circle (.132); 
\draw [thick] (-.384,1.067) circle (.132); 
\draw [thick] (.384,1.067) circle (.132); 
\draw [thick] (1.115,-0.201) circle (.132); 
\draw [thick] (-1.115,-0.201) circle (.132); 
\end{tikzpicture}
\caption{Left: tree-like structure of the barycentric chart.  Right: An analogous circle packing of a triangle.}
\label{fig:tree}
\end{figure}
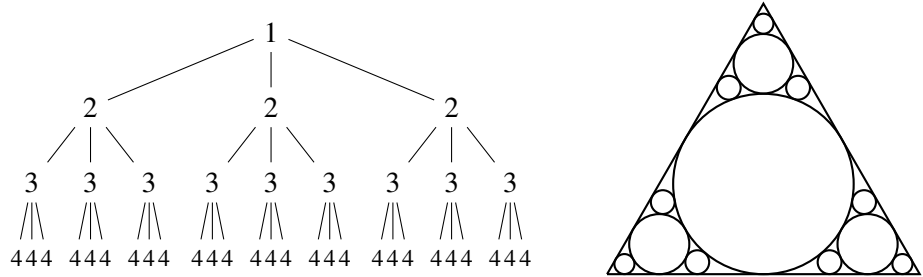

~
The rational points of the original spiderweb remain rational in barycentric coordinates.
 E.g., the tangency points between ellipses in the main chain are now
$$
\left(\!\left(\;  
\frac{1}{4n^2+1},\,  
\frac{2n(n-1)}{4n^2+1},\,  
\frac{2n(n+1)}{4n^2+1}
\; \right)\!\right)
$$
(and the permutations).  
Points of tangency between the region of depth 1 and ellipses in the chain become:
$$
\left(\!\left(\;  
\frac{1}{2(n^2-n+1)},\,  
\frac{(n-1)^2}{2(n^2-n+1)},\,  
\frac{n^2}{2(n^2-n+1)}
\; \right)\!\right)
$$
Finally, the points of tangency between the diagonal ellipses are
$$
\left(\!\left(\;  
\frac{n(n-1)}{n^2-n+1},\,  
\frac{1/2}{n^2-n+1},\,  
\frac{1/2}{n^2-n+1}
\; \right)\!\right)
$$

\newpage

\section{Additional bits}

\begin{figure}[H]
\centering
\includegraphics[scale=.38]{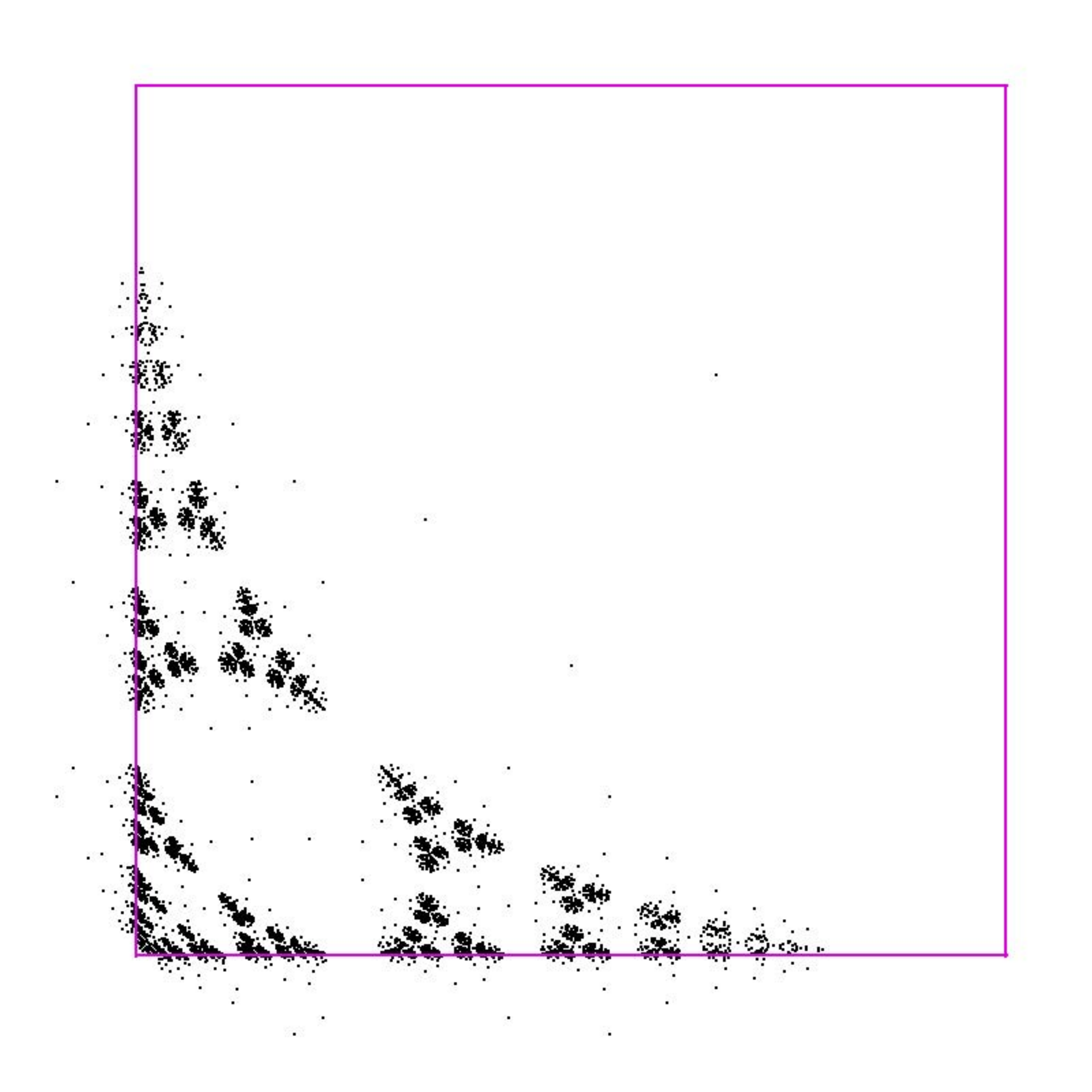} 
\caption{\small The dust of the tricycles (seeds) the Apollonian Window plotted in the chart  (Dust2)}
\label{fig:dust}
\end{figure}

\paragraph{A. Apollonian disk packings in the chart.} 
Let $\hat p$ denote the Apollonian circle packing (up to similarity) 
determined by tricycle $p\in \mathcal T$.
Define the equivalence relation
$$
p\sim q \qquad \hbox{if}\qquad \hat p = \hat q
$$
The quotient $\mathcal T/\sim$ is the moduli space of Apollonian disk packings.
Each individual Apollonian disk packing understood
as the set of all tricycles it contains
may be drawn as a dust of points in $\mathcal T$.
Figure \ref{fig:dust} shows the equivalence class of tricycles corresponding to the Apollonian Window.
Since every point of $\mathcal T$ leads to an Apollonian packing, 
the space splits into an uncountable number of sets of countable many points.

Among problems that are interesting and easy to state but not necessarily simple Is:
Is there a continuous path in $\mathcal T$ that includes all 
Apollonian packings without repetitions?

\begin{figure}[H]
\centering
\includegraphics[scale=.24]{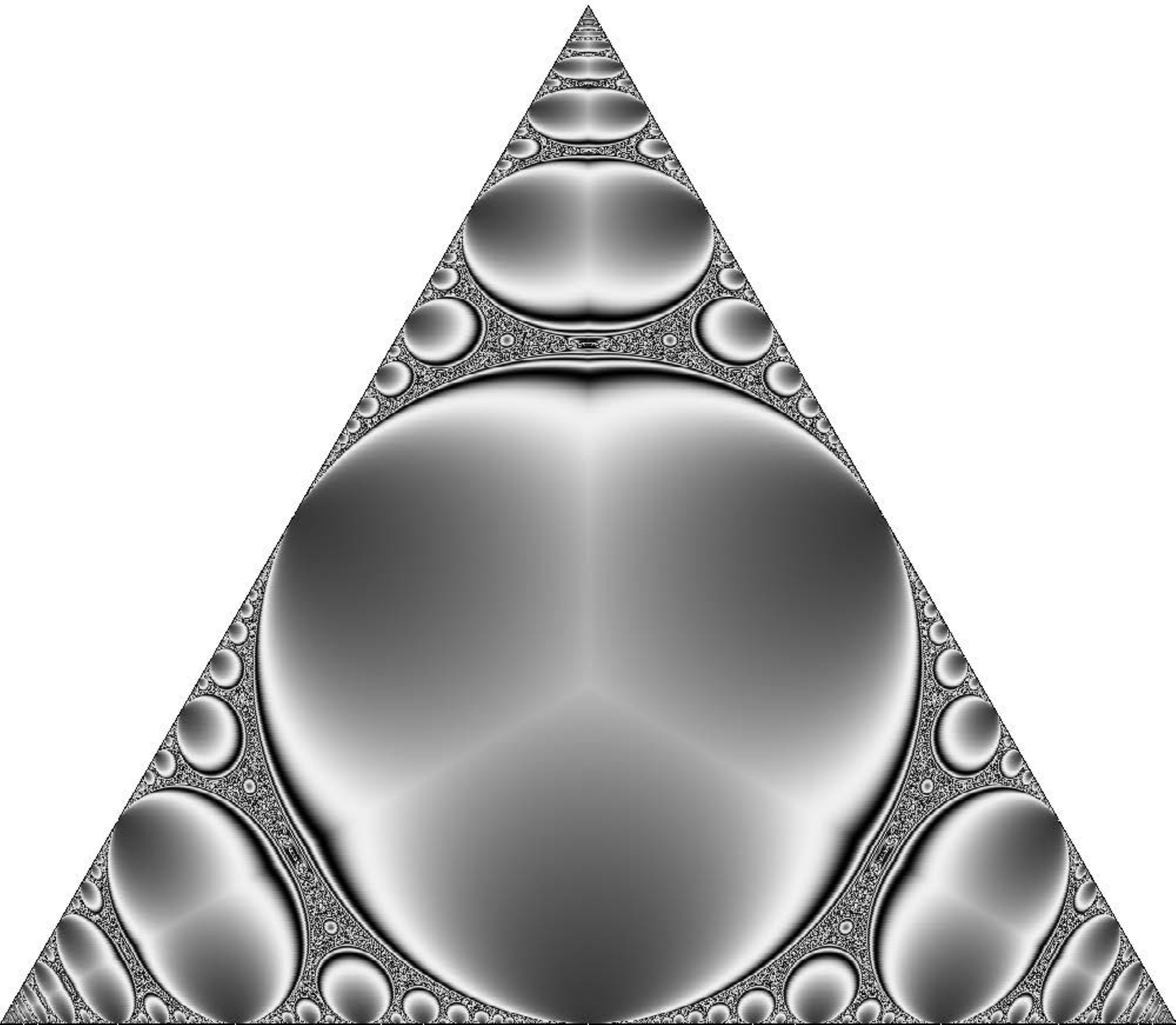} 
\caption{\small Size depth value, sinusoidal change}
\label{fig:T-size-BW}
\end{figure}

\paragraph{B. Size function.}  A scalar function related to that of depth is a function 
$$
f: \mathcal T \ \to \ \mathbb R \cup\{\infty\}
$$
associating to every tricycle
the radius (or curvature) of the size of the greatest circle of the Apollonian disk packing it determines.
A version for the barycentric setup of the chart is shown in Figure \ref{fig:T-size-BW}.
 It actually shows the composition $g=\sin\circ f$. The fine pattern in some regions 
 results from the interaction between resolution of the drawing with the resolution of the 
 rapidly changing values of $g$.

\paragraph{C. Squaring}

\begin{figure}[H]
\centering
\includegraphics[scale=.21]{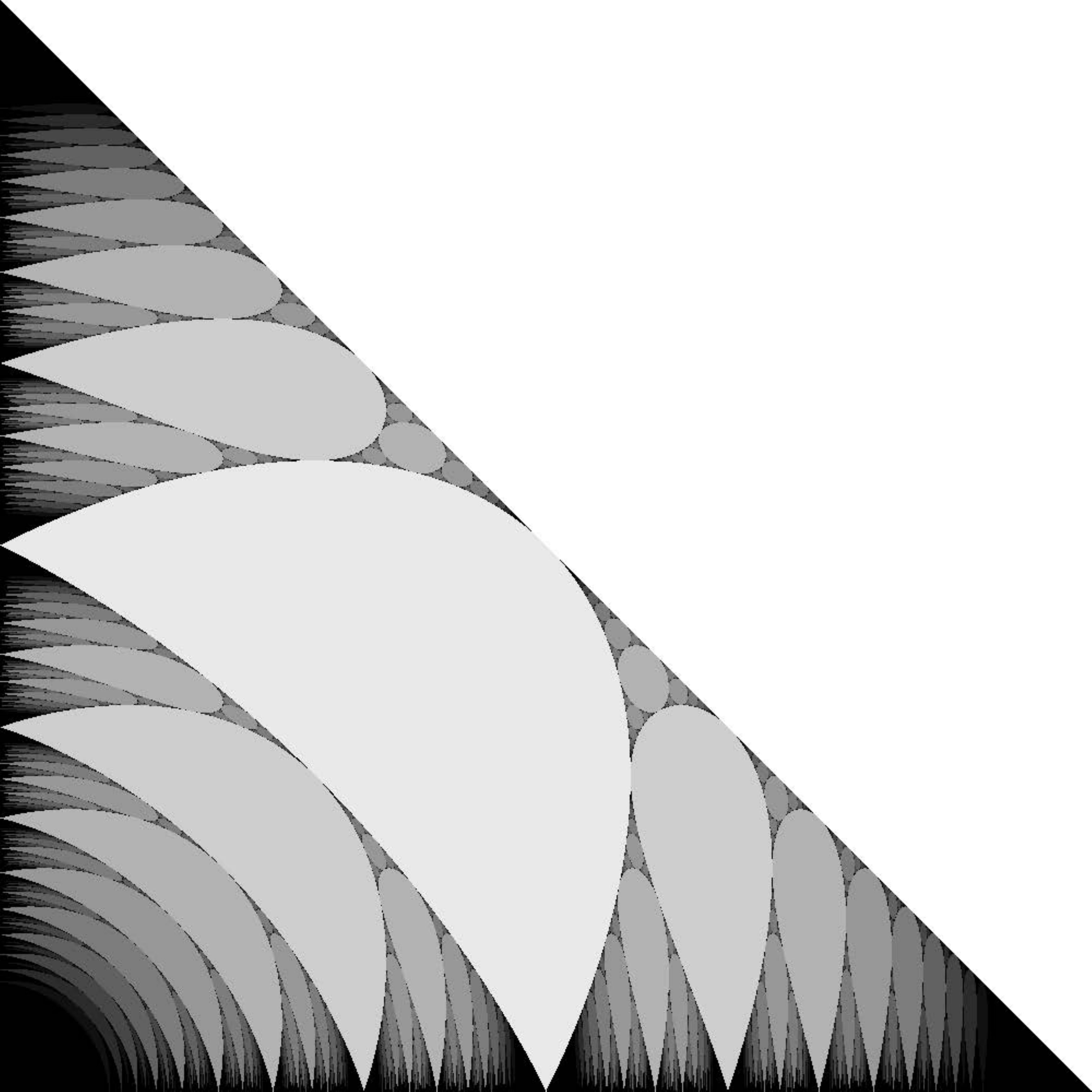} 
\caption{\small Quadratic deformation}
\label{fig:drzewa}
\end{figure}

Since the points of tangency are all squares of rational numbers, one might think that re-scaling the figure 
so that the squares are brought to non-squares will ``straighten'' the figure and the ellipses will become circles.
Modifying the code and rerunning the program 
invalidates the guess.
But a dramatic image it produces is shown in Figure \ref{fig:drzewa}. 

~\\
{\bf E. Open questions and challenges.}  
\\\\
{\bf 1.}  Find the exact equations of all ellipses in the fractal and a consistent way to present them
as function of the ellipses addresses (see Figure \ref{fig:Sierpinski1}). 
\\[7pt]
{\bf 2.}  The rationality of tangency points 
prompt further investigations into plausible number-theoretic relations.
\\[7pt]
{\bf 3.}  Is there a single (conformal?) map that will transform all the ellipses of the spiderweb into circles?

~


\section*{Appendix A: The algorithm}

The pseudo-code of the algorithm is presented below:

\begin{figure}[H]
\hrule

\vspace{.14in}
{\small
\begin{verbatim} 
for(int n = 0;  n < 1000; n++) 
  { for(int m=0; m <= n; m++) 
      { new = 1;   depth=0; 
        T[0] = n; T[1] = m, T[2] = 1000; 

        while ( (new > 0) && (depth < 21) )
           { depth=depth+1;   
              T = sort(T);
              float a = T[0];  float b = T[1]; float c = T[2]; 
              new = a + b + c - 2*sqrt(a*b + b*c + c*a);
              T[2] = new;
	           }  //end while 

        colorRGB( 30*depth, 30*depth, 30*depth ); 
        draw point(n,m);  draw point(m,n);

       }  // end of ``for m''
   }      // end of ``for n''
\end{verbatim}
}

\hrule

~

\caption{The code for drawing the chart.}
\label{fig:code}
\end{figure}



\section*{Appendix B: More ellipse equations}

$$
\begin{array}{llll}
1L_{1,1}: \qquad&  \left( x+y \right)^2 +1    &\quad = \quad &4xy + 2x+2y \\[10pt]
2L_{1,1}: &  \left( 2x+2y \right)^2+\dfrac{1}{4}    &\quad = \quad  &4xy + 2x+2y \\[10pt]
3L_{1,1}: &   \left( 3x+3y\right)^2  + \dfrac{1}{9}  &\quad = \quad  &4xy + 2x+2y \\[10pt]
3L_{2,1}: &  \left( \dfrac{3}{2} x +2y \right)^2   +\dfrac{4}{9}  &\quad = \quad  &4xy + 2x +  \dfrac{4}{3} y \\[10pt]
4L_{1,1}: &  \left( 4x +4y \right)^2 +\dfrac{1}{16}  &\quad = \quad  &4xy + 2x+2y  \\[10pt]
4L_{3,2}: &  \left(  \dfrac{14}{9}x + \dfrac{11}{6} y \right)^2  + \dfrac{4}{9}    
                                &\quad = \quad  &4xy +  \dfrac{52}{27} x + \dfrac{14}{9} y  \\[10pt]
4L_{4,1}: &  \left( \dfrac{4}{3} x+  2y \right)^2    + \dfrac{9}{16}    &\quad = \quad &4xy +  2x+y  \\[10pt]
4L_{3,1}:  &  \left( \dfrac{5}{3}x + \dfrac{11}{5} y \right)^2 + \dfrac{9}{25} 
                              &\quad = \quad  &4xy +   2x +  \dfrac{34}{25} y \\[10pt]
4L_{2,1}: &  \left( \dfrac{5}{2}x+\dfrac{16}{5}y \right)^2 +\dfrac{4}{25}    &\quad = \quad  &4xy +  2x+ \dfrac{44}{25}y\\[10pt]
%
%
5L_{1,1}: &   (5x+5y)^2+\frac{1}{25}   &\quad = \quad  & 4xy +2x+2y   \\[10pt]
5L_{8,1}: &    (\frac{5}{4} x+2y)^2  +\frac{16}{25}  
       &\quad = \quad  & 4xy +2x + \frac{4}{5} y   \\[10pt]
5L_{5,1}:  &    (\frac{7}{4} x+\frac{16}{7} y)^2+\frac{16}{49}  
        &\quad = \quad  & 4xy+2x+ \frac{68}{49} y   \\[10pt]
5L_{5,4}:  &    (\frac{19}{12} x + \frac{16}{9} y)^2 +\frac{4}{9} 
         &\quad = \quad  & 4xy  + \frac{17}{9} x+ \frac{44}{27} y   \\[10pt]
5L_{4,1}:  &    (\frac{7}{3} x + \frac{19}{7} y)^2  +\frac{9}{49}   
           &\quad = \quad  & 4xy  +2x+ \frac{82}{49} y  \\[10pt]
5L_{6,1}:  &    (\frac{8}{5} x+ \frac{11}{5} y)^2    + \frac{25}{64} 
           &\quad = \quad  & 4xy +2x+ \frac{5}{4} y  \\[10pt]
5L_{6,3}:  &    (\frac{3}{2} x+ \frac{11}{6} y)^2  + \frac{121}{256} 
         &\quad = \quad  &  4xy  + \frac{31}{16} x + \frac{71}{48} y   \\[10pt]
5L_{3,1}:  &    (\frac{8}{3} x+3y)^2+\frac{9}{64}  
            &\quad = \quad  &  4xy  +2x + \frac{7}{4} y  \\[10pt]
5L_{3,1}:   &  ( \frac{7}{5} x+\frac{73}{35} y)^2 + \frac{25}{49}    
            &\quad = \quad  &   4xy + 2x+ \frac{50}{49} y  \\[10pt]
5L_{2,1}:  &    (\frac{7}{2} x+\frac{26}{7} y)^2+\frac{4}{49}  
           &\quad = \quad  &  4xy+2x+ + \frac{92}{49} y  \\[10pt]
5L_{7,2}:  &   ( \frac{26}{19} x + \frac{73}{38} y)^2 +\frac{196}{361}  
           &\quad = \quad  &  4xy  + \frac{716}{361} x + \frac{422}{361} y  \\[10pt]
5L_{5,3}:  &    (\frac{27}{17} x+\frac{97}{51} y)^2+\frac{121}{289}  
            &\quad = \quad  & 4xy + \frac{562}{289} x+ \frac{1334}{867} y  \\[10pt]
5L_{3,2}:  &    (\frac{38}{15} x + \frac{27}{10} y)^2   +\frac{4}{25}  
               &\quad = \quad  &  4xy  + \frac{148}{75} x+ \frac{46}{25} y  \\[10pt]
5L_{5,2}:   &    ( \frac{38}{23} x+\frac{97}{46})^2 +\frac{196}{529}  
              &\quad = \quad  &   4xy + \frac{1052}{529} x+ \frac{758}{529} y  \\[10pt]
\end{array}
$$

%

%
%

%
%
%
%
%

\newpage

\section*{Appendix C: A remark on the ``Apollonian''}

\begin{figure}[h]
\centering
\begin{tikzpicture}[scale=1.32]
\draw [thick] (-1,-1) circle (1);
\draw [thick] (-1,1) circle (1);
\draw [thick] (1,-1) circle (1);
\draw [thick] (1,1) circle (1);
\draw [thick] (0,0) circle (.412);
\draw [] (0,.519) circle (.1055);
\draw [] (.519,0) circle (.1055);
\draw [] (0,-.519) circle (.1055);
\draw [] (-.519,0) circle (.1055);
\draw [] (0,.675) circle (.048);
\draw [] (0,-.675) circle (.048);
\draw [] (.675,0) circle (.048);
\draw [] (-.675,0) circle (.048);
\draw [] (.1215,.438) circle (.03884);
\draw [] (-.1215,.438) circle (.03884);
\draw [] (.1215,-.438) circle (.03884);
\draw [] (-.1215,-.438) circle (.03884);
\draw [] (.438,.1215) circle (.03884);
\draw [] (.438,-.1215) circle (.03884);
\draw [] (-.438,.1215) circle (.03884);
\draw [] (-.438,-.1215) circle (.03884);
\draw [] (0,.752) circle (.027);
\draw [] (0,-.752) circle (.027);
\draw [] (.752,0) circle (.027);
\draw [] (-.752,0) circle (.027);
\node at (-1, 1) [scale=1, color=black] {\sf a};
\node at (1, 1) [scale=1, color=black] {\sf b};
\node at (0, 0) [scale=1, color=black] {\sf c};
\end{tikzpicture}
\qquad
\begin{tikzpicture}[scale=16, rotate=-0]%
\draw [thick] (3/6, 4/6) circle (1/6);
\draw [thick] (6/11,8/11) circle (1/11);
\draw [thick] (6/15,10/15) circle (1/15);
\draw [thick] (7/14,8/14) circle (1/14);
\draw [thick] (14/23, 14/23) circle (1/23);
\draw [thick] (11/26, 20/26) circle (1/26);
\draw [thick] (14/35, 20/35) circle (1/35);
\draw [thick] (27/42, 28/42) circle (1/42);
\draw [thick] (22/47, 38/47) circle (1/47);
\draw [thick] (30/51, 28/51) circle (1/51);
\draw [thick] (22/59, 44/59) circle (1/59);
\draw [thick] (30/71, 38/71) circle (1/71);
\draw [thick] (46/71, 50/71) circle (1/71);
\draw [thick] (27/74, 44/74) circle (1/74);
\draw [thick] (41/86, 56/86) circle (1/86);
\draw [thick] (54/95, 50/95) circle (1/95);
\draw [thick] (49/110, 80/110) circle (1/110);
\draw [thick] (57/102, 64/102) circle (1/102);
\draw [thick] (70/107, 68/107) circle (1/107);
\draw [thick] (75/122, 68/122) circle (1/122);
\draw [thick] (57/134, 80/134) circle (1/134);
\draw [thick] (97/158, 104/158) circle (1/158);
\draw [thick] (105/182, 104/182) circle (1/182);
\draw [thick] (99/210, 140/210) circle (1/210);
\draw [thick] (102/219, 140/219) circle (1/219);
\end{tikzpicture}
\caption{Left: A non-Apollonian packing.
Right: another Apollonian packing.}
\label{fig:square}
\end{figure}
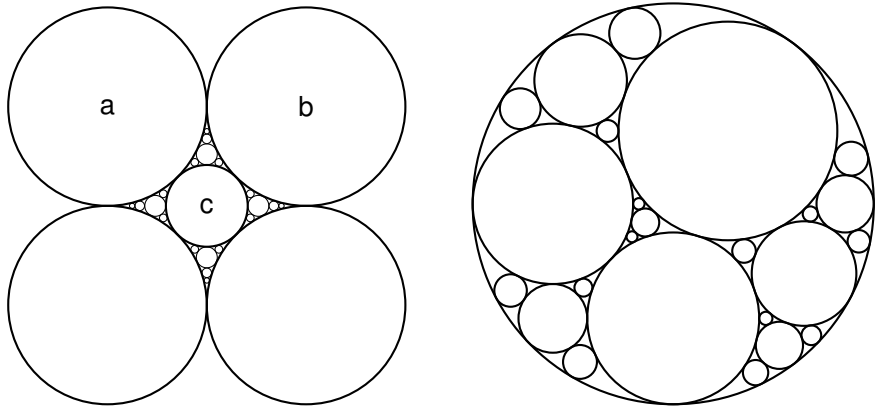

Figure \ref{fig:square} shows a generic Apollonian disk packing (left) and a non-Apollonian disk packing
(only one of the two completions of the triple (a,b,c) belongs to the packing.)

\noindent
An Apollonian disk packing is a collection of disks  (some of possibly of non-positive curvature)
such that the following {\bf Apollonian rules} are satisfied:  \
(1)  It contains a tricycle;  (2)  No disks overlap; \
(3)  For any three mutually tangent circles, both Descartes solutions also belong to  the packing.

\section*{Acknowledgments}

I am indebted to Philip Feinsilver for his comments on this manuscript.


\newpage
\begin{center}
{\bf \Huge Additional images}
\end{center}

~

Some color images that are too ink-draining to be part of the main body of the text.
Putting them here allows one to prevent printing them.

~

\begin{figure}[H]
\includegraphics[scale=.37]{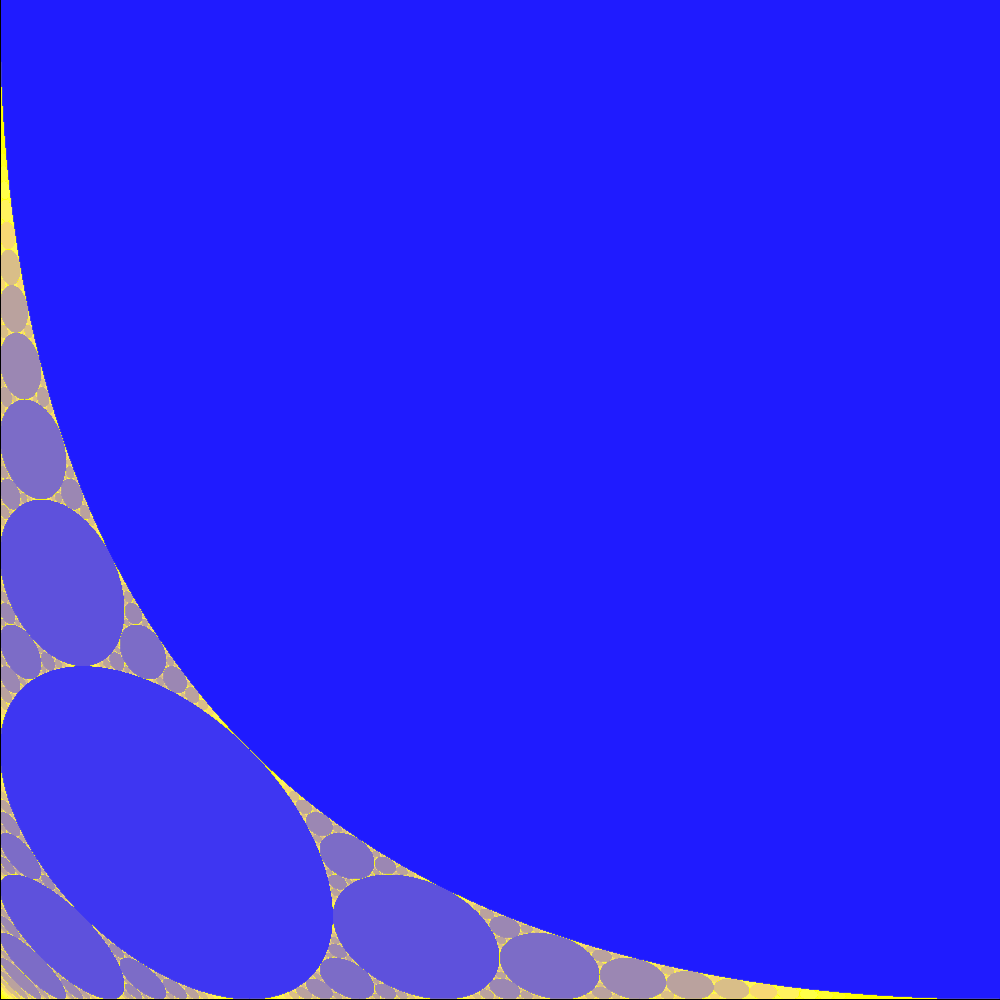}
\caption{The bluer the smaller the value of the depth function.}
\label{fig:end1}
\end{figure}

\begin{figure}[H]
\centering
\includegraphics[scale=.49]{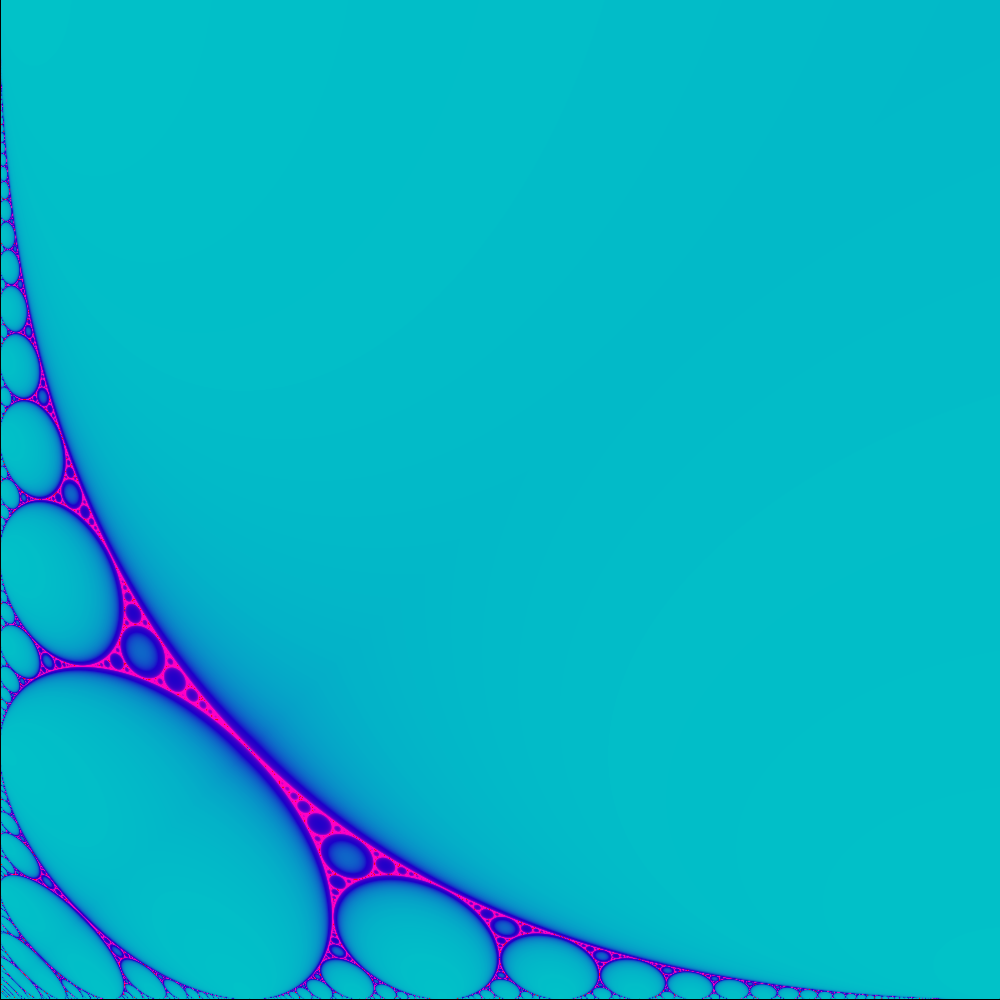}
\caption{The Apollonian size function. The purple color marks the higher values.}
\label{fig:end2}
\end{figure}

\begin{figure}[H]
\includegraphics[scale=.37]{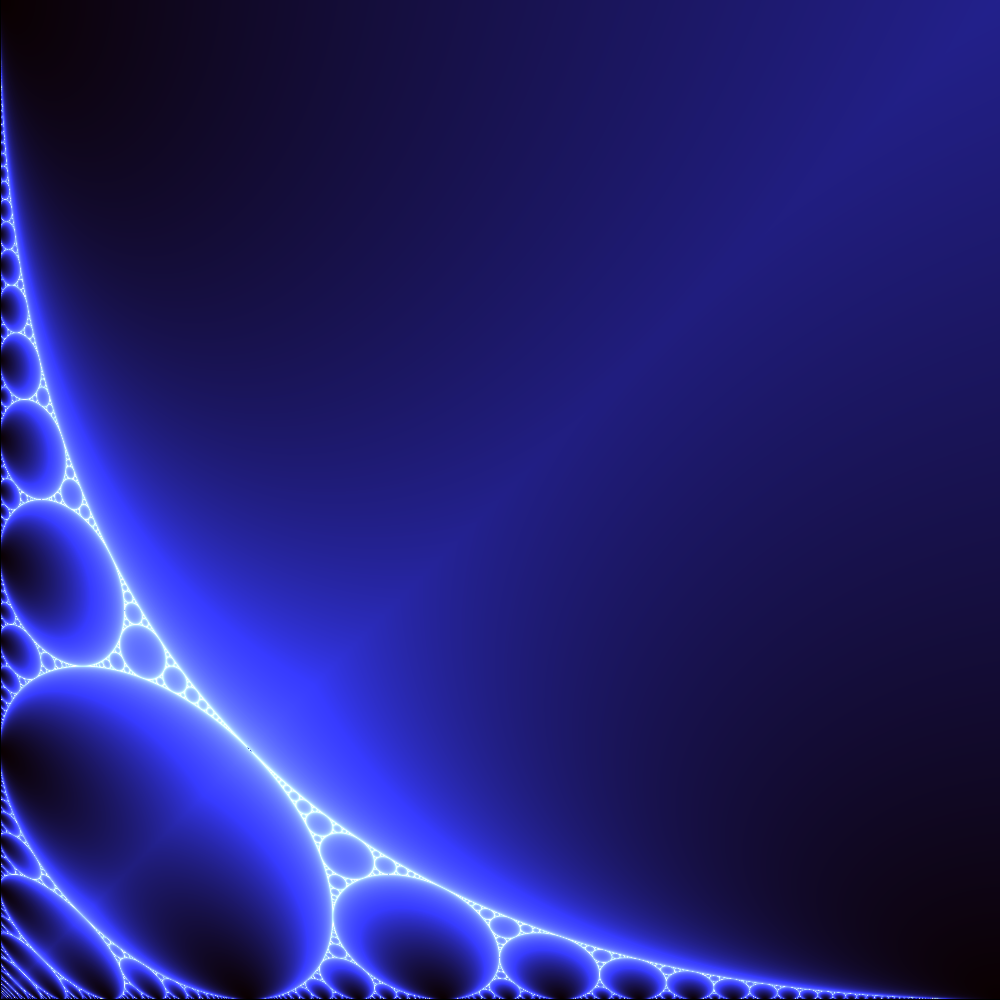}
\caption{The Apollonian size function again. Lighter pixels for the greater values.}
\label{fig:end3}
\end{figure}

\begin{figure}[H]
\includegraphics[scale=.37]{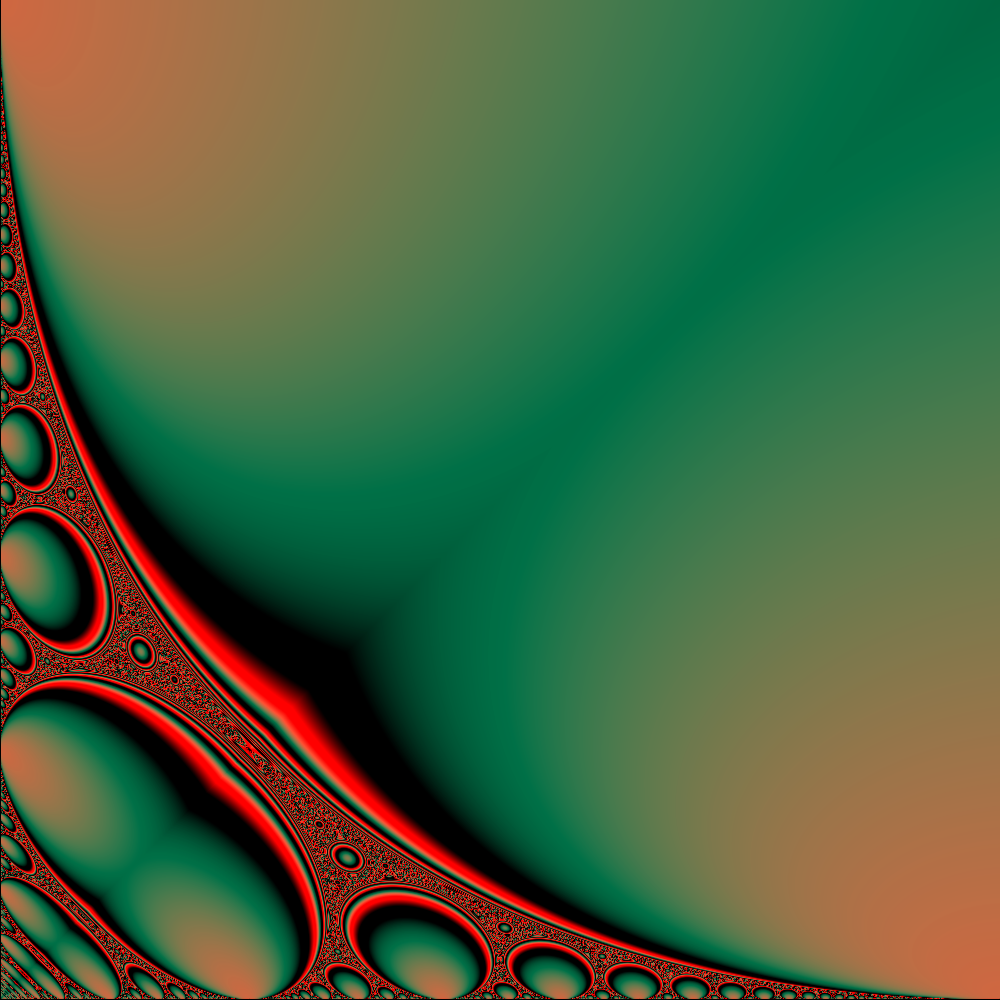}
\caption{The Apollonian size function represented via swiftly changing colors in a cyclic way.}
\label{fig:end4}
\end{figure}

\end{document}